\documentclass[a4paper,12pt]{article}     

\usepackage{graphicx}		  % to deal with graphicals
\usepackage{amsmath,amssymb}	% to deal with mathematics

\newtheorem{lemma}{Lemma}[section]
\newtheorem{theorem}[lemma]{Theorem}
\newtheorem{proposition}[lemma]{Proposition}
\newtheorem{corollary}[lemma]{Corollary}
\newtheorem{assumption}{Assumption}

\newcommand\LABEL[1]{\label{#1}}

\def\authorfont{\footnotesize}

\def\ccode#1{\par
\vspace*{8pt}
{\authorfont{\leftskip18pt\rightskip\leftskip
\noindent #1\par}}\par}

\newenvironment{Proof}{
\hspace*{-9mm}
{ \it Proof.}}
{\hfill {$\square$}\vspace{1.5em}}

\begin{document}

\begin{center}{
{\Large 
 Properties of minimal charts and
 their applications VI: 
 the graph $\Gamma_{m+1}$ in a chart $\Gamma$ of type $(m;2,3,2)$}
\vspace{10pt}
\\ 
Teruo NAGASE and Akiko SHIMA\footnote{The second author is supported by JSPS KAKENHI Grant Number 18K03309.}
}
\end{center}

%%%%%%%%%%%%%%
%%%%%%%%%%%%%%%
%%%%%%%%%%%%%%%%% 
%%%%%%%%%%%%%%%%
\date{2020.03.20}

%%%%%%%%%%%%%%%%%%%%%%%%%%%%%%%
%%%%   abstract   %%%%%%%%%%%%%
%%%%%%%%%%%%%%%%%%%%%%%%%%%%%%%
\begin{abstract}
Let $\Gamma$ be a chart,
and we denote by $\Gamma_m$
the union of all the edges of label $m$.
A chart $\Gamma$ is of type $(m;2,3,2)$
if $w(\Gamma)=7$,
$w(\Gamma_m\cap\Gamma_{m+1})=2$,
$w(\Gamma_{m+1}\cap\Gamma_{m+2})=3$,
and $w(\Gamma_{m+2}\cap\Gamma_{m+3})=2$
where 
$w(G)$ is the number of white vertices in $G$.
In this paper, we prove that if there is 
a minimal chart $\Gamma$ of 
type $(m;2,3,2)$,
then
each of $\Gamma_{m+1}$ and $\Gamma_{m+2}$
contains one of three kinds of graphs.
In the next paper,
we shall prove
that there is no minimal chart of type $(m;2,3,2)$.
\end{abstract}

%
%
%%%%   Make Content   %%%%%%%%%%
%%%%%%%%%%%%%%%%%%%%%%%%%%%%%%%%
%\tableofcontents
%
%
%%%%%%%%%%%%%%%%%%%%%%%%%%%%
%%%%%%%%%%%%%%%%%%%%%%%%%%%%

\ccode{2010 Mathematics Subject Classification. Primary 57Q45; Secondary 57Q35.}
\ccode{ {\it Key Words and Phrases}. surface link, chart, white vertex. }

\setcounter{section}{0}
\section{Introduction}

%\large
%\baselineskip 20pt

Charts are oriented labeled graphs in a disk (see  \cite{KnottedSurfaces},\cite{BraidBook}, and see Section~\ref{s:Prel}  for the precise definition of charts).
From a chart, we can construct an oriented closed surface 
embedded in 4-space ${\Bbb R}^4$ 
 (see \cite[Chapter 14, Chapter 18 and Chapter 23]{BraidBook}). 
A C-move 
is a local modification between two charts
in a disk (see Section~\ref{s:Prel} for C-moves).
A C-move between two charts induces 
an ambient isotopy between oriented closed surfaces 
corresponding to the two charts.

We will work in the PL category or smooth category. All submanifolds are assumed to be locally flat.
In \cite{ONS},
we showed that there is no minimal chart with exactly five vertices
 (see Section~\ref{s:Prel} for the precise definition of minimal charts). 
Hasegawa proved that there exists a minimal chart with exactly
six white vertices \cite{H1}. 
This chart represents a 2-twist spun trefoil.
In \cite{INS} and \cite{NST},
we investigated minimal charts with exactly four white vertices.
In this paper, 
we investigate properties of minimal charts and
need to prove that
there is no minimal chart with exactly seven white vertices
(see \cite{ChartApp1},\cite{ChartAppII},
\cite{ChartAppIII},\cite{ChartAppIV},
\cite{ChartAppV}, \cite{ChartAppVII},
\cite{ChartAppVIII}).

Let $\Gamma$ be a chart.
For each label $m$, we denote by $\Gamma_m$
the union of all the edges of label $m$.

Now we define a type of a chart:
Let $\Gamma$ be a chart with at least one white vertex, 
and $n_1,n_2,\dots,n_k$ integers.
The chart $\Gamma$ is of {\it type $(n_1,n_2,\dots,n_k)$} if there exists a label $m$ of $\Gamma$ satisfying the following three conditions:
\begin{enumerate}
\item[(i)] For each $i=1,2,\dots, k$, 
the chart $\Gamma$ contains exactly $n_{i}$ white vertices in $\Gamma_{m+i-1}\cap \Gamma_{m+i}$.
\item[(ii)] If $i<0$ or $i>k$, then $\Gamma_{m+i}$ does not contain any white vertices.
\item[(iii)] Both of the two subgraphs $\Gamma_m$ and $\Gamma_{m+k}$ contain at least one white vertex.
\end{enumerate}
If we want to emphasize the label $m$,
then we say that $\Gamma$ is of {\it type $(m;n_1,n_2,\dots,n_k)$}. 
Note that $n_1\ge1$ and $n_k\ge1$ by the condition (iii).

We proved in \cite[Theorem 1.1]{ChartAppII} that
if there exists a minimal $n$-chart $\Gamma$ with exactly seven white vertices,
then $\Gamma$ is a chart of 
type $(7),(5,2),(4,3),(3,2,2)$ or $(2,3,2)$ 
(if necessary we change the label
$i$ by $n-i$ for all label $i$).
In \cite{ChartAppV},
we showed that
there is no minimal chart of type $(3,2,2)$.
In this paper and \cite{ChartAppVII},
we shall show 
the following.

\begin{theorem}
\LABEL{MainTheorem} 
{\rm (\cite[Theorem 1.1]{ChartAppVII})}
There is 
no minimal chart of 
type $(2,3,2)$.
\end{theorem}

In the future paper \cite{ChartAppVIII},
we shall show there is no minimal chart of type
$(7),(5,2),(4,3)$.
Therefore we shall show that
there is no minimal chart with exactly seven white vertices.

An edge in a chart is called 
a {\it terminal edge}
if it has
a white vertex and a black vertex.

 In our argument  we often construct a chart $\Gamma$. 
On the construction of a chart $\Gamma$, for a white vertex $w\in\Gamma_m$ for some label $m$,  
among the three edges of $\Gamma_m$ 
containing $w$, 
if one of the three edges is a terminal edge 
(see Fig.~\ref{Fig01}(a) and (b)), 
then we remove the terminal edge and
put a black dot at the center of the white vertex  as shown in Fig.~\ref{Fig01}(c).
Namely
Fig.~\ref{Fig01}(c) means 
Fig.~\ref{Fig01}(a) or 
Fig.~\ref{Fig01}(b).
We call the vertex in Fig.~\ref{Fig01}(c) 
a {\it BW-vertex} with respect to $\Gamma_m$.

%%%%%%%%%%%%%%%%%%
%%%%%%%%%%%%%%%%%% Figure
%%%%%%%%%%%%%%%%%%
\begin{figure}[hbt]
\centerline{\includegraphics{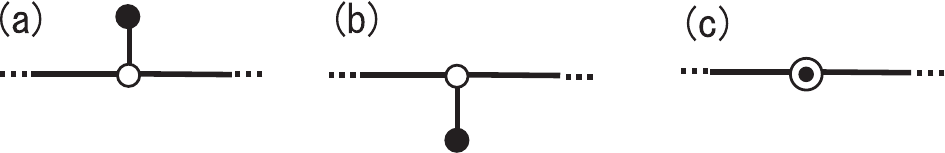}}
\caption{\LABEL{Fig01}
(a),(b) White vertices in terminal edges.
(c) A BW-vertex.}
\end{figure}

For example, the graph as shown in Fig.~\ref{Fig02}(a) means one of the four graphs as shown in
Fig.~\ref{Fig02}(b),(c),(d),(e).

%%%%%%%%%%%%%%%%%%
\begin{figure}[hbt]
\centerline{\includegraphics{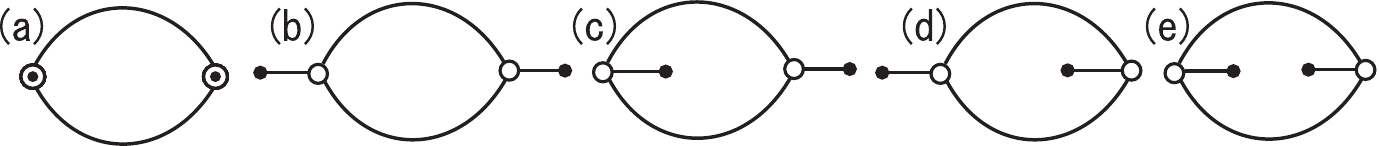}}
\caption{\LABEL{Fig02}
Graphs with two white vertices.}
\end{figure}

The three graphs in Fig.~\ref{Fig03}
are examples of graphs in $\Gamma_m$ for a chart $\Gamma$
and a label $m$.
We call 
a {\it $\theta$-curve},
an {\it oval},
a {\it skew $\theta$-curve} the three graphs as shown
in Fig.~\ref{Fig03}(a),(b),(c)
respectively.

%%%%%%%%%%%%%%%%%%
%%%%%%%%%%%%%%%%%% Figure
%%%%%%%%%%%%%%%%%%
\begin{figure}[htb]
\centerline{\includegraphics{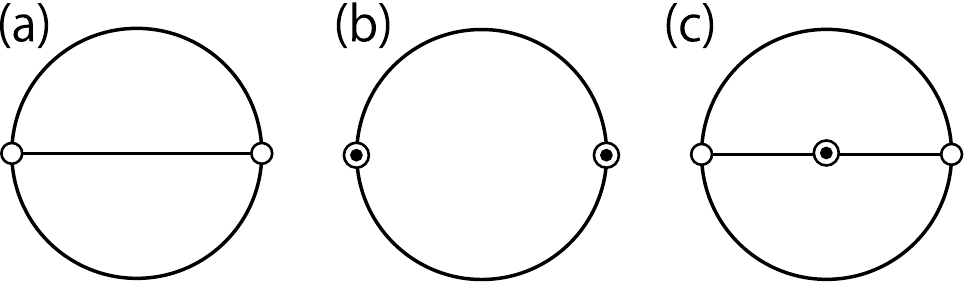}}
\caption{\LABEL{Fig03}
(a) A $\theta$-curve. (b) An oval. 
(c) A skew $\theta$-curve.}
\end{figure}

Let $X$ be a set in a chart $\Gamma$.
Let
 $$w(X)=\text{the number of white vertices in $X$.}$$

Let $\Gamma$ be a chart of type $(m;2,3,2)$.
Then $w(\Gamma)=7,w(\Gamma_m\cap\Gamma_{m+1})=2,
w(\Gamma_{m+1}\cap\Gamma_{m+2})=3,
w(\Gamma_{m+2}\cap\Gamma_{m+3})=2$.
Thus $w(\Gamma_{m+1})=5$ and $w(\Gamma_{m+2})=5$.
First we shall show the following lemma.

\begin{lemma}
\LABEL{GammaM+1}
Let $\Gamma$ be a minimal chart of type $(m;2,3,2)$.
Then each of $\Gamma_{m+1}$ and $\Gamma_{m+2}$
contains one of nine graphs as shown 
in  Fig.~\ref{Fig04}, or
the union of 
a $\theta$-curve and a skew $\theta$-curve, or
the union of 
an oval and a skew $\theta$-curve.
\end{lemma}

%%%%%%%%%%%%%%%%%%
%%%%%%%%%%%%%%%%%% Figure
%%%%%%%%%%%%%%%%%%
\begin{figure}[htb]
\centerline{\includegraphics{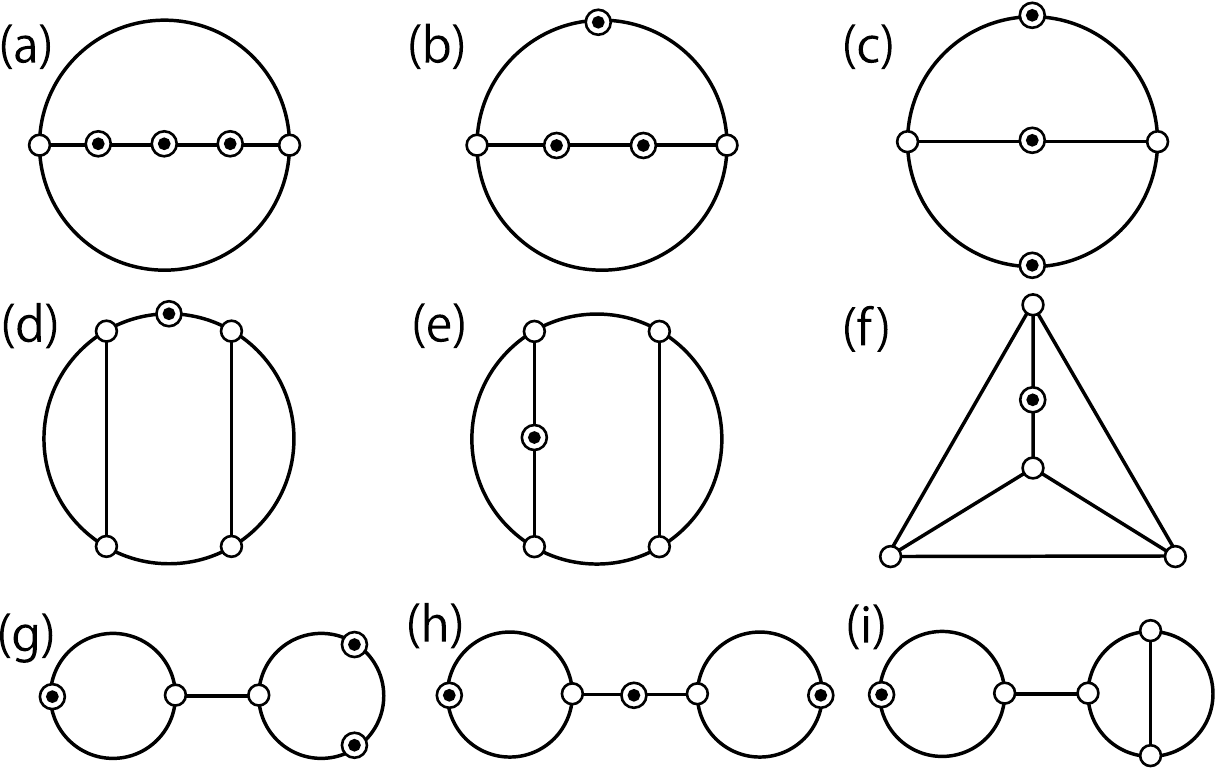}}
\caption{\LABEL{Fig04}
(a),(b),(c) Graphs with three black vertices.
(d),(e),(f) Graphs with one black vertex.
(g),(h) Graphs with three black vertices.
(i) A graph with one black vertex.}
\end{figure}

In this paper,
the following is the main result.

\begin{theorem}
\LABEL{OvalTypeGTypeH}
If there exists a minimal chart $\Gamma$ of 
type $(m;2,3,2)$,
then each of $\Gamma_{m+1}$ and $\Gamma_{m+2}$
contains either
the union of an oval and a skew $\theta$-curve, 
or 
 one of two graphs as shown 
in  Fig.~\ref{Fig04}$($g$)$,$($h$)$.
\end{theorem}

%%%%%%%%%%%%%%%%%%%%%%%%%
%%%%%%%%%%%%%%%%%%%%%%%%%%
%%%%%%%%%%%%%%%%%%%%%%%%%%

%\newpage

The paper is organized as follows.
In Section~\ref{s:Prel},
we define charts and minimal charts.
In Section~\ref{s:ConnectedComponent},
we investigate connected components of $\Gamma_m$
with five white vertices
for a minimal chart $\Gamma$.
We shall show Lemma~\ref{GammaM+1}.
In Section~\ref{s:ThetaCurve},
we shall show that
neither $\Gamma_m$ nor $\Gamma_{m+3}$
contains a $\theta$-curve
for any minimal chart $\Gamma$ of type $(m;2,3,2)$
(i.e. both of $\Gamma_m$ and $\Gamma_{m+3}$
contain ovals) (see Corollary~\ref{NoThetaGammaM}).
In Section~\ref{s:Oval}.
we investigate an oval of label $m$
for a minimal chart $\Gamma$.
In Section~\ref{s:WhiteGammaM+1GammaM+2},
we investigate white vertices 
in an oval of label $m$ 
for a minimal chart $\Gamma$ of type $(m;2,3,2)$.
In Section~\ref{s:GammaM+1GammaM+2},
we shall show that
for any minimal chart $\Gamma$ of type $(m;2,3,2)$,
the graph $\Gamma_{m+1}$ contains none of
the five graphs as shown in 
Fig.~\ref{Fig04}(a),(d),(e),(f),(i),
and neither does $\Gamma_{m+2}$.
Moreover
we shall show that
neither $\Gamma_{m+1}$ nor $\Gamma_{m+2}$
contains a $\theta$-curve.
In Section~\ref{s:ROfamilies},
we consider a minimal chart $\Gamma$ of type $(m;2,3,2)$
such that $\Gamma_{m+1}$ contains either an oval, or 
one of the four graphs as shown in
Fig.~\ref{Fig04}(b),(c),(g),(h).
We investigate 
that the chart $\Gamma$ contains 
what kind of pseudo charts.
In Section~\ref{s:IOC},
we shall show that
neither $\Gamma_{m+1}$ nor $\Gamma_{m+2}$
contains the graph as shown in Fig.~\ref{Fig04}(c)
for any minimal chart $\Gamma$ of type $(m;2,3,2)$.
In Section~\ref{s:ShiftingLemma},
we shall show that
neither $\Gamma_{m+1}$ nor $\Gamma_{m+2}$
contains the graph as shown in Fig.~\ref{Fig04}(b)
for any minimal chart $\Gamma$ of type $(m;2,3,2)$.
We obtain the main theorem (Theorem~\ref{OvalTypeGTypeH}).

%%%%%%%%%%%%%%%%%%%%%%%%%%%%%%%%%%%%%%
%%%%%%%%%%%%%%%%%%%%%%%%%%%%%%%%%%%%%%
%%%%%%%%%%%%%%%%%%%%%%%%%%%%%%%%%%%%%%
%%%%%%%%%%%%%%%%%%%%%%%%%%%%%%%%%%%%%%
%%%%%%%%%%%%%%%%%%%%%%%%%%%%%%%%%%%%%%
%%%%%%%%%%%%%%%%%%%%%%%%%%%%%%%%%%%%%%

%\newpage
\section{Preliminaries}
\LABEL{s:Prel}

In this section, 
we introduce 
the definition of charts and its related words.

Let $n$ be a positive integer.
An $n$-{\it chart}  
(a braid chart of degree $n$ \cite{KnottedSurfaces}
or a surface braid chart of degree $n$ \cite{BraidBook}) 
is 
an oriented labeled graph in the interior of a disk,
which may be empty 
or
have closed edges without vertices
satisfying the following four conditions
(see Fig.~\ref{Fig05}):
\begin{enumerate}
\item[(i)] 
Every vertex has degree $1$, $4$, or $6$.
\item[(ii)] 
The labels of edges are 
in $\{1,2,\dots,n-1\}$.
\item[(iii)]
In a small neighborhood of
each vertex of degree $6$,
there are six short arcs,
three consecutive arcs are
oriented inward 
and
the other three are outward,
and
these six are labeled $i$ and $i+1$
alternately for some $i$,
where the orientation and label of
each arc are inherited from
the edge containing the arc.
\item[(iv)]
For each vertex of degree $4$,
diagonal edges have the same label
and
are oriented coherently,
and the labels $i$ and $j$ of
the diagonals satisfy $|i-j|>1$.
\end{enumerate}
We call a vertex of degree $1$ a {\it black vertex},
a vertex of degree $4$ a {\it crossing}, and 
a vertex of degree $6$ a {\it white vertex}
respectively.

Among six short arcs
in a small neighborhood of
a white vertex,
a central arc of each three consecutive arcs
oriented inward (resp. outward) 
is called a   
{\it middle arc} at the white vertex
(see Fig.~\ref{Fig05}(c)).
For each white vertex $v$, 
there are two middle arcs at $v$ 
in a small neighborhood of $v$.
An edge is said to be {\it middle at} a white vertex $v$ if it contains a middle arc at $v$.

Let $e$ be an edge connecting $v_1$ and $v_2$.
If $e$ is oriented from $v_1$ to $v_2$,
then we say that 
$e$ is oriented {\it outward at $v_1$}
and {\it inward at $v_2$}
%%%%%%%%%%%%%%%%%%

%%%%%%%%%%%%%%%%%%

\begin{figure}[htb]
\begin{center}
\includegraphics{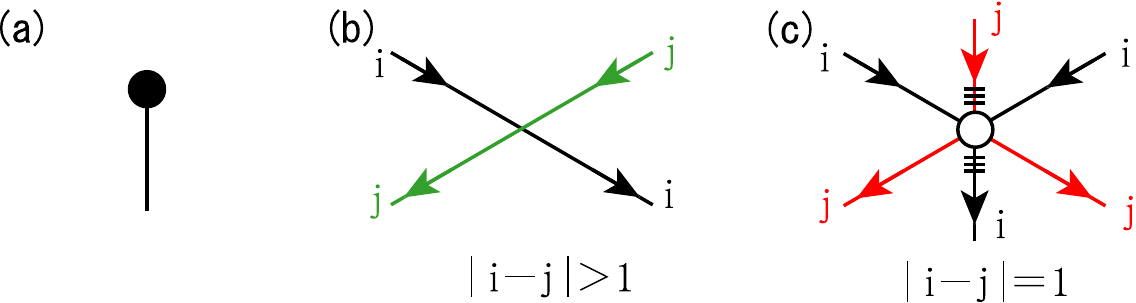}
\end{center}
\caption{ \LABEL{Fig05} (a) A black vertex. (b) A crossing. (c) A white vertex. 
Each arc with three transversal short arcs is a middle arc at the white vertex. }
\end{figure}

Now {\it C-moves} are local modifications 
of charts as shown in Fig.~\ref{Fig06}
(cf. \cite{KnottedSurfaces}, 
\cite{BraidBook} and \cite{Tanaka}).
Two charts are said to be {\it C-move equivalent}  if there exists
a finite sequence of C-moves 
which modifies one of the two charts 
to the other.

%%%%%%%%%%%%%%%%%%%
\begin{figure}
\begin{center}
\includegraphics{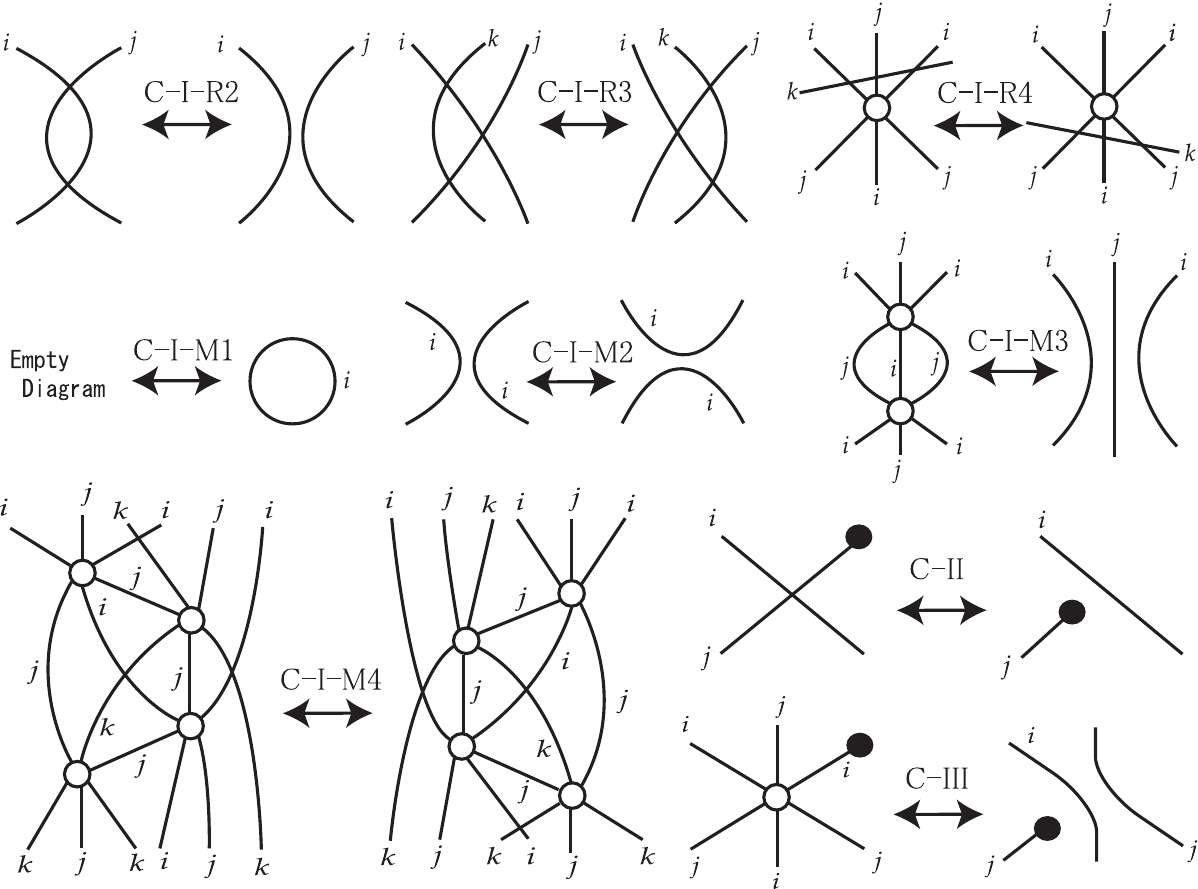}
\end{center}
\caption{ \LABEL{Fig06} For the C-III move, 
the edge with the black vertex is not middle at
a white vertex in the left figure. }
\end{figure}
%%%%%%%%%%%%%%%%%%

An edge in a chart is called 
a {\it free edge}
if it has
two black vertices.

For each chart $\Gamma$,
let $w(\Gamma)$ and $f(\Gamma)$ be the number of white vertices, and the number of free edges respectively.
The pair $(w(\Gamma), -f(\Gamma))$ is called a {\it complexity} of the chart (see \cite{BraidThree}).
A chart $\Gamma$ is called a {\it minimal chart} if its complexity is minimal among the charts C-move equivalent to the chart $\Gamma$ with respect to the lexicographic order of pairs of integers.

We showed the difference of a chart in a disk and in a 2-sphere (see \cite[Lemma 2.1]{ChartApp1}).
This lemma follows from that there exists a natural one-to-one correspondence between $\{$charts in $S^2\}/$C-moves and $\{$charts in $D^2\}/$C-moves, conjugations
(\cite[Chapter 23 and Chapter 25]{BraidBook}).
To make the argument simple, we assume that 
the charts lie on the 2-sphere instead of the disk.
\begin{assumption}
In this paper,
all charts are contained in the $2$-sphere $S^2$.
\end{assumption}
We have the special point in the 2-sphere $S^2$, called the point at infinity,
 denoted by $\infty$.
In this paper, all charts are contained in a disk such that the disk 
does not contain the point at infinity $\infty$.

Let $\Gamma$ be a chart,
and $m$ a label of $\Gamma$. 
A {\it hoop} is a closed edge of $\Gamma$ without vertices 
(hence without crossings, neither).
A {\it ring} is a simple closed curve in $\Gamma_m$ containing a crossing but not containing any white vertices.
A hoop is said to be {\it simple} 
if one of the two complementary domains
of the hoop
does not contain any white vertices.

We can assume that
all minimal charts $\Gamma$
satisfy the following four conditions 
(see \cite{ChartApp1},\cite{ChartAppII},\cite{ChartAppIII},
\cite{StI}):

%%%%%%%%%%%%%%%%%%%%%%%%%%%%%%%%%%%%%%%%%
\begin{assumption}
\LABEL{AssumeTerminal}
If an edge of $\Gamma$
contains a black vertex,
then the edge is a free edge 
or a terminal edge.
Moreover 
any terminal edge contains a middle arc.
\end{assumption}

%%%%%%%%%%%%%%%%%%%%%%%%%%%%%%%%%%%%%%%
\begin{assumption}
\LABEL{NoSimpleHoop}
All free edges and simple hoops in $\Gamma$ 
are moved into a small neighborhood $U_\infty$ 
of the point at infinity $\infty$. 
Hence
we assume that 
$\Gamma$ does not contain free edges
nor simple hoops, 
otherwise mentioned. 
\end{assumption}

%%%%%%%%%%%%%%%%%%%%%%%%%%%%%%%%%%%%%%%%	
\begin{assumption}
\LABEL{Ring}
Each complementary domain of
any ring and hoop must contain 
at least one white vertex. 
\end{assumption}

\begin{assumption}
\LABEL{Infinity}
The point at infinity $\infty$ is moved in any complementary domain of $\Gamma$.
\end{assumption}

In this paper
for a set $X$ in a space
we denote 
the interior of $X$,
the boundary of $X$ and
the closure of $X$
by Int$X$, $\partial X$
and $Cl(X)$
respectively.

%%%%%%%%%%%%%%%%%%%%
%%%%%%%%%%%%%%%%%%%%%
%%%%%%%%%%%%%%%%%%%%%
%\newpage
\section{Connected components of $\Gamma_m$}
\LABEL{s:ConnectedComponent}
In this section,
we investigate connected components of $\Gamma_m$
with five white vertices
for a minimal chart $\Gamma$.
We shall show Lemma~\ref{GammaM+1}.

\begin{lemma}
\LABEL{OriBWvertex}
{\rm (\cite[Lemma 3.1]{ChartAppV})}
In a minimal chart $\Gamma$,
for each BW-vertex in $\Gamma_m$,
the two edges of label $m$ containing the BW-vertex
are oriented inward or outward at the BW-vertex
simultaneously
if each of the two edges is not a terminal edge
$($see Fig.~\ref{Fig07}$)$.
\end{lemma}

%%%%%%%%%%%%%%%%%%
%%%%%%%%%%%%%%%%%% Figure
%%%%%%%%%%%%%%%%%%
\begin{figure}[hbt]
\centerline{\includegraphics{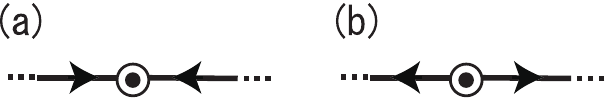}}
\caption{\LABEL{Fig07}
BW-vertices.}
\end{figure}

Let $\Gamma$ be a chart,
and $m$ a label of $\Gamma$. 
A {\it loop} is a simple closed curve in $\Gamma_m$ with exactly one white vertex
(possibly with crossings).

\begin{lemma}
\LABEL{LemmaWithTerminal3}
{\rm (\cite[Lemma 3.2]{ChartAppV})}
Let $\Gamma$ be a minimal chart,
and $m$ a label of $\Gamma$.
Let $G$ be a connected component of $\Gamma_m$.
Then we have the following.
\begin{enumerate}
\item[{\rm (a)}] If $1\le w(G)$, then $2\le w(G)$.
\item[{\rm (b)}] If $1\le w(G)\le 3$
and $G$ does not contain any loop, 
then $G$ is one of three graphs as shown 
in Fig.~\ref{Fig03}.
\end{enumerate}
\end{lemma}

The following lemma is easily shown.
Thus we omit the proof.

\begin{lemma}
\LABEL{LemmaWithLoop}
Let $G$ be a $3$-regular graph in $S^2$.
Then we have the following.
\begin{enumerate}
\item[{\rm (a)}] The graph $G$ contains exactly 
an even number of vertices. 
\item[{\rm (b)}]
If $G$ has at most four vertices, then 
$G$ is one of seven graphs as shown 
in Fig.~\ref{Fig03}$($a$)$ 
and Fig.~\ref{Fig08}.
\end{enumerate}
\end{lemma}

%%%%%%%%%%%%%%%%%%
%%%%%%%%%%%%%%%%%% Figure
%%%%%%%%%%%%%%%%%%
\begin{figure}[htb]
\centerline{\includegraphics{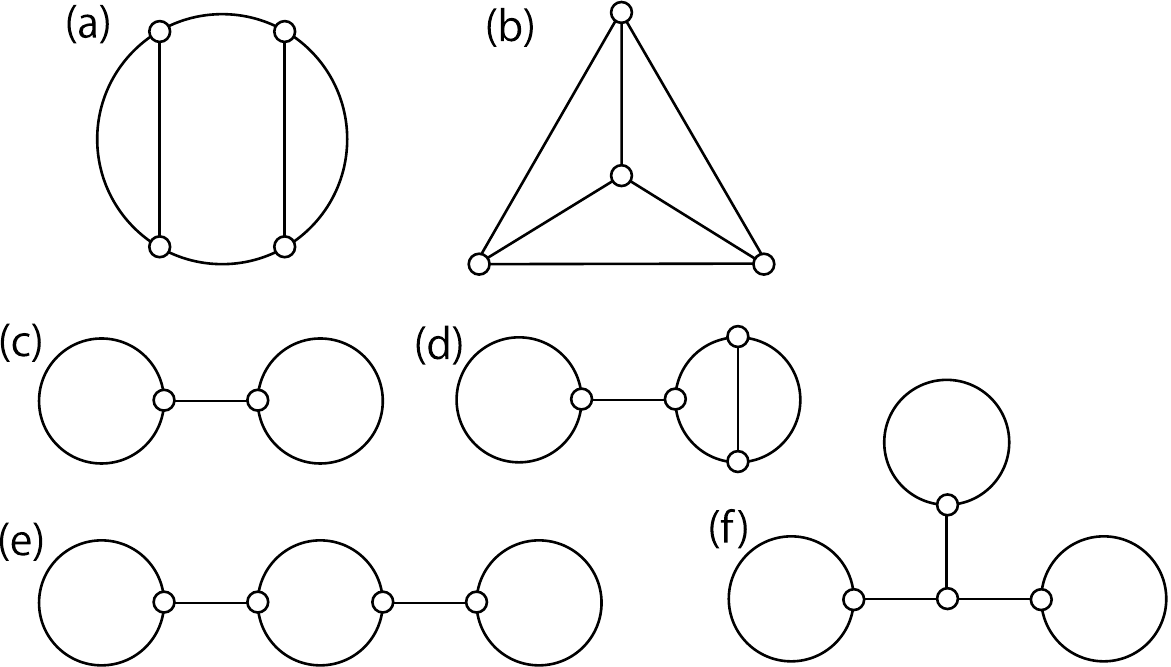}}
\caption{\LABEL{Fig08}
(a),(b) Graphs without loops.
(c),(d),(e),(f) Graphs with loops.}
\end{figure}

\begin{lemma}
\LABEL{LemmaWithTerminal}
Let $\Gamma$ be a minimal chart,
and $m$ a label of $\Gamma$.
Let $G$ be a connected component of $\Gamma_m$.
If $w(G)=5$ and $G$ has no loop,
then $G$ is one of nine graphs as shown 
in  Fig.~\ref{Fig04}.
\end{lemma}

\begin{Proof}
By Assumption~\ref{AssumeTerminal},
each terminal edge is middle at a white vertex.
Thus each white vertex in $\Gamma_m$ is 
contained in at most one terminal edge of label $m$.
Hence 
\begin{enumerate}
\item[$(1)$]the graph $G$ is obtained 
from a simple closed curve or 
a 3-regular graph (possibly with loops) 
by adding BW-vertices.
\end{enumerate}

Now we shall show that 
$G$ contains at least one black vertex.
If not,
then the graph $G$ is a 3-regular graph on $S^2$.
By Lemma~\ref{LemmaWithLoop}(a),
the graph $G$ contains exactly 
an even number of white vertices.
This contradicts the fact $w(G)=5$.
Hence $G$ contains at least one black vertex.
Thus
\begin{enumerate}
\item[$(2)$]
$G$ contains at least one BW-vertex.
\end{enumerate}

{\bf Claim.}
The graph $G$ is obtained from a 3-regular graph by adding BW-vertices.

Suppose that all white vertices in $G$ are BW-vertices.
Then the graph $G$ is obtained from a simple closed curve 
by adding BW-vertices.
By Lemma~\ref{OriBWvertex},
in a minimal chart,
for each BW-vertex in $\Gamma_m$,
the two edges of label $m$ containing the BW-vertex
are oriented inward or outward at the BW-vertex
simultaneously
if each of the two edges is not a terminal edge.
Hence the orientation of edges must change at BW-vertices.
Thus $G$ contains exactly an even number of BW-vertices.
This contradicts the fact $w(G)=5$.
Hence $G$ contains a white vertex not a BW-vertex.
Thus by (1),
Claim holds. {\hfill {$\square$}\vspace{1.5em}}

By $w(G)=5$, (2) and Claim,
the graph $G$ is obtained by adding BW-vertices
from a 3-regular graph with at most four vertices.
Hence by Lemma~\ref{LemmaWithLoop}(b),
the graph $G$ is obtained by adding BW-vertices
from one of seven graphs as shown 
in Fig.~\ref{Fig03}(a) 
and Fig.~\ref{Fig08}.
Since $G$ has no loop with $w(G)=5$,
the graph $G$ is not
obtained from the graphs 
as shown in Fig.~\ref{Fig08}(e),(f).

Now the graph $G$ is on the 2-sphere $S^2$.
Hence if $G$ is obtained from the graph as shown 
in Fig.~\ref{Fig03}(a), 
then 
the graph $G$ is one of three graphs as shown
in Fig.~\ref{Fig04}(a),(b),(c).
If $G$ is obtained from the graph 
as shown in Fig.~\ref{Fig08}(a), 
then 
the graph $G$ is one of two graphs as shown 
in Fig.~\ref{Fig04}(d),(e).
If $G$ is obtained from the graph 
as shown in Fig.~\ref{Fig08}(b), 
then 
the graph $G$ is the graph as shown  
in Fig.~\ref{Fig04}(f).
If $G$ is obtained from 
 the graph as shown in Fig.~\ref{Fig08}(c),
then 
the graph $G$ is one of two graphs as shown  
in Fig.~\ref{Fig04}(g),(h).
If $G$ is obtained from the graph as shown in 
Fig.~\ref{Fig08}(d),
then 
the graph $G$ is the graph as shown 
in Fig.~\ref{Fig04}(i).
Therefore 
$G$ is one of nine graphs as shown  
in Fig.~\ref{Fig04}.
We complete the proof of Lemma~\ref{LemmaWithTerminal}.
\end{Proof}

\begin{lemma}
\LABEL{GammaFiveWhite}
Let $\Gamma$ be a minimal chart,
and $m$ a label of $\Gamma$.
If $w(\Gamma_m)=5$ and $\Gamma_m$ has no loop,
then $\Gamma_m$ contains one of the following graphs:
\begin{enumerate}
\item[{\rm (a)}] one of nine graphs as shown 
in  Fig.~\ref{Fig04}, or
\item[{\rm (b)}] the union of 
a $\theta$-curve and a skew $\theta$-curve, or
\item[{\rm (c)}] the union of 
an oval and a skew $\theta$-curve.
\end{enumerate}
\end{lemma}

\begin{Proof}
First we shall show that  there exist at most two 
connected components of $\Gamma_m$ with
white vertices.
Suppose that there exist at least three 
connected components 
$G_1,G_2,G_3$ of $\Gamma_m$ with
$w(G_i)\ge1$ for each $i=1,2,3$.
Then by Lemma~\ref{LemmaWithTerminal3}(a),
we have $w(G_i)\ge2$ for each $i=1,2,3$.
Thus
$$5=w(\Gamma_m)\ge w(G_1)+w(G_2)+w(G_3)\ge 2+2+2=6.$$
This is a contradiction.
Hence there exist at most two 
connected components of $\Gamma_m$ with
white vertices.

Suppose that
there exists a connected component $G_1$ of $\Gamma_m$
with $w(G_1)=5$.
Since $\Gamma_m$ has no loop,
by Lemma~\ref{LemmaWithTerminal}
the graph $G_1$ is one of nine graphs as shown 
in  Fig.~\ref{Fig04}.

Suppose that
there exists two connected components $G_1,G_2$ 
of $\Gamma_m$ with $w(G_1)\ge1,w(G_2)\ge1$ and
$w(G_1)+w(G_2)=5$.
Then by Lemma~\ref{LemmaWithTerminal3}(a),
we have $w(G_1)\ge2,w(G_2)\ge2$.

Without loss of generality
we can assume $2\le w(G_1)\le w(G_2)$.
Since $w(G_1)+w(G_2)=5$,
we have $w(G_1)=2$ and $w(G_2)=3$.
Since $\Gamma_m$ has no loop, 
by Lemma~\ref{LemmaWithTerminal3}(b)
the graph $G_1$ is a $\theta$-curve or an oval, and
the graph $G_2$ is a skew $\theta$-curve.
\end{Proof}

\begin{lemma}$(${\rm \cite[Theorem 1.1]{ChartAppIV}}$)$
\LABEL{LemmaNoLoop}
There is no loop in any minimal chart with exactly seven white vertices.
\end{lemma}

By Lemma~\ref{GammaFiveWhite} and
Lemma~\ref{LemmaNoLoop},
we have Lemma~\ref{GammaM+1}.

%%%%%%%%%%%%%%%%%%%%%%%%%%%%
%%%%%%%%%%%%%%%%%%%%%%%%%%%%
%%%%%%%%%%%%%%%%%%%%%%%%%%%%

%\newpage
\section{$\theta$-curves}
\LABEL{s:ThetaCurve}

In this section
we shall show that
neither $\Gamma_m$ nor $\Gamma_{m+3}$
contains a $\theta$-curve
for any minimal chart $\Gamma$ of type $(m;2,3,2)$
(i.e. both of $\Gamma_m$ and $\Gamma_{m+3}$
contain ovals) (see Corollary~\ref{NoThetaGammaM}).

Let $\Gamma$ be a chart, 
and $m$ a label of $\Gamma$.
Let $L$ be the closure of a connected component 
of the set obtained by taking out 
all the white vertices from $\Gamma_m$.
If $L$ contains at least one white vertex
but does not contain any black vertex,
then $L$ is called an {\it internal edge of label $m$}.
Note that an internal edge may contain a crossing of $\Gamma$.

Let $\Gamma$ be a chart. 
Let $D$ be a disk 
such that 
\begin{enumerate}
\item[(1)] the boundary $\partial D$ consists of an internal edge $e_1$ of label $m$ and an internal edge $e_2$ of label ${m+1}$, and 
\item[(2)] any edge containing a white vertex in $e_1$ does not intersect the open disk Int$D$.
\end{enumerate}
Note that $\partial D$ may contain crossings.
Let $w_1$ and $w_2$ be the white vertices in $e_1$. 
If the disk $D$ satisfies one of the following conditions, then $D$ is called  {\it a lens of type $(m,m+1)$}
(see Fig.~\ref{Fig09}):
\begin{enumerate}
	\item[(i)] Neither $e_1$ nor $e_2$ contains a middle arc. 
	\item[(ii)] One of the two edges $e_1$ and $e_2$ contains middle arcs at both white vertices $w_1$ and $w_2$ simultaneously.
\end{enumerate}

%%%%%%%%%%%%%%%%%% Figure
%%%%%%%%%%%%%%%%%%
\begin{figure}[htb]
\centerline{\includegraphics{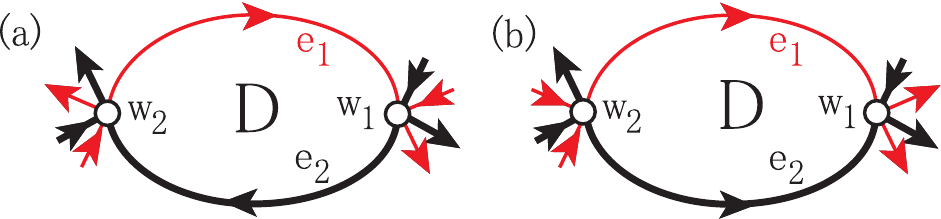}}
\caption{\LABEL{Fig09}
Lenses.}
\end{figure}

\begin{lemma}{\rm (\cite[Corollary 1.3]{ChartAppII})}
\LABEL{NoLens}
 There is no lens in any minimal chart with 
at most seven white vertices.
\end{lemma}

\begin{lemma}
\LABEL{ThetaGammaM}
Let $\Gamma$ be a minimal chart,
and $m$ a label of $\Gamma$.
Suppose that $\Gamma_m$ contains
a $\theta$-curve $G$.
If $\Gamma$ has no lens,
and if the two white vertices in $G$ are
contained in $\Gamma_{m+\varepsilon}$ for some
$\varepsilon\in\{+1,-1\}$,
then $w(\Gamma_{m+\varepsilon})\ge6$.
\end{lemma}

\begin{Proof}
Let $w_1,w_2$ be the white vertices in $G$,
and $e$ the internal edge of label $m$ in $G$ 
middle at $w_1$.
Without loss of generality we can assume that
\begin{enumerate}
\item[$(1)$]
the edge $e$ is oriented from $w_1$ to $w_2$.
\end{enumerate}
Then the other two internal edges in $G$
are oriented from $w_2$ to $w_1$
(see Fig.~\ref{Fig10}(a)).
Thus 
\begin{enumerate}
\item[$(2)$]
the edge $e$ is middle at $w_1,w_2$.
\end{enumerate}

The $\theta$-curve $G$ divides $S^2$ into three disks.
Let $D_1,D_2$ be two of the three disks
with $D_1\cap D_2=\partial D_1\cap \partial D_2=e$.
Let $e_1,e_2$ be internal edges
(possibly terminal edges) of label $m+\varepsilon$ 
at $w_1$ in $D_1,D_2$ respectively.
Since $e$ is middle  at $w_1$ by (2),
neither $e_1$ nor $e_2$ is middle at $w_1$.
By Assumption~\ref{AssumeTerminal},
neither $e_1$ nor $e_2$ is a terminal edge.
Hence
$e_1$ and $e_2$ contain white vertices 
different from $w_1$, say $w_3,w_4$.

We shall show $w_3\not=w_2$ and $w_4\not=w_2$.
If $w_3=w_2$,
then the edge $e_1$ separates the disk $D_1$ into 
two disks.
One of the two disks contains the edge $e$.
By (2),
the disk is a lens.
This contradicts the condition that
$\Gamma$ has no lens.
Hence $w_3\not=w_2$.
Similarly we can show $w_4\not=w_2$.

Let $e_1',e_2'$ be internal edges 
(possibly terminal edges) of label $m+\varepsilon$ at $w_2$
in $D_1,D_2$ respectively.
By using (2),
we can show similarly that
$e_1',e_2'$ contain white vertices different from $w_2$,
say $w_3',w_4'$.

We shall show that
$w(\Gamma_{m+\varepsilon}\cap{\rm Int}D_1)\ge2$.
There are two cases:
$w_3\not=w_3'$ and $w_3=w_3'$.

If $w_3\not=w_3'$,
then  $w(\Gamma_{m+\varepsilon}\cap{\rm Int}D_1)\ge2$.

Suppose $w_3=w_3'$ (see Fig.~\ref{Fig10}(b)).
Let $e_1''$ be an internal edge
(possibly a terminal edge) of label $m+\varepsilon$
at $w_3$ 
different from $e_1,e_1'$.
By (1) and (2),
the edge $e_1$ is oriented from $w_1$ to $w_3$
and 
the edge $e_1'$ is oriented from $w_3$ to $w_2$.
Thus $e_1''$ is not middle at $w_3$.
Hence
by Assumption~\ref{AssumeTerminal},
the edge $e_1''$ is not a terminal edge.
Thus the edge $e_1''$ contains a white vertex 
different from $w_3$.
Thus $w(\Gamma_{m+\varepsilon}\cap{\rm Int}D_1)\ge2$.

Similarly we can show 
 $w(\Gamma_{m+\varepsilon}\cap{\rm Int}D_2)\ge2$.
Finally we have
$$w(\Gamma_{m+\varepsilon})\ge 
w(\Gamma_{m+\varepsilon}\cap G)+
w(\Gamma_{m+\varepsilon}\cap {\rm Int}D_1)+
w(\Gamma_{m+\varepsilon}\cap {\rm Int}D_2)
\ge 2+2+2=6$$
\end{Proof}

%%%%%%%%%%%%%%%%%%
%%%%%%%%%%%%%%%%%% Figure
%%%%%%%%%%%%%%%%%%
\begin{figure}[hbt]
\centerline{\includegraphics{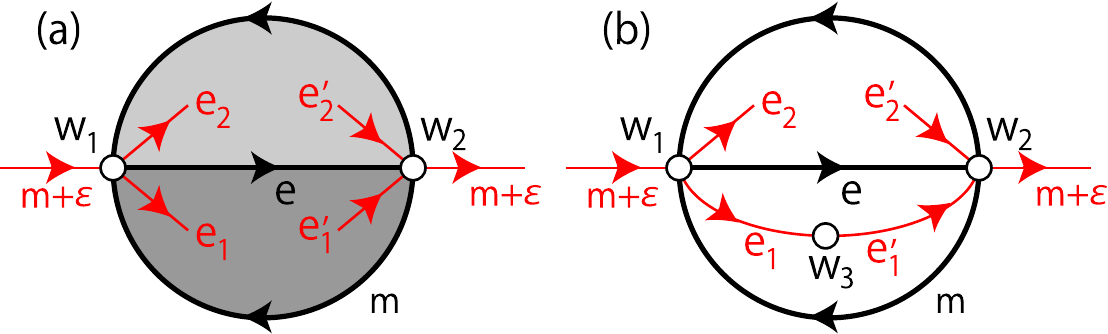}}
\caption{\LABEL{Fig10}
(a) The dark gray region is the disk $D_1$,
the light gray region is the disk $D_2$.
(b) Both of $e_1$ and $e_1'$ contain 
the white vertex $w_3$.
}
\end{figure}

\begin{corollary}
\LABEL{NoThetaGammaM}
Let $\Gamma$ be a minimal chart of type $(m;2,3,2)$.
Then both of  $\Gamma_{m}$ and $\Gamma_{m+3}$
contain ovals.
\end{corollary}

\begin{Proof}
Since $\Gamma$ is of type $(m;2,3,2)$,
we have $w(\Gamma_m)=2$.
Since the graph $\Gamma_m$ does not contain any loop
by Lemma~\ref{LemmaNoLoop},
Lemma~\ref{LemmaWithTerminal3}(b) implies that
the graph $\Gamma_m$
contains one of the two graphs
as shown in Fig.~\ref{Fig03}(a) and (b).
Hence the graph $\Gamma_m$
contains a $\theta$-curve or an oval.
If the graph $\Gamma_m$ contains a $\theta$-curve,
then we have $w(\Gamma_{m+1})\ge6$ 
by Lemma~\ref{ThetaGammaM}.

On the other hand, 
since $\Gamma$ is of type $(m;2,3,2)$,
we have $w(\Gamma_{m+1})=5$.
This is a contradiction.
Thus the graph $\Gamma_m$
contains an oval.

Similarly
we can show that
 the graph $\Gamma_{m+3}$
contains an oval.
\end{Proof}

%%%%%%%%%%%%%%%%%%%%%%%%%%%%
%%%%%%%%%%%%%%%%%%%%%%%%%%%%
%%%%%%%%%%%%%%%%%%%%%%%%%%%%

%\newpage
\section{Ovals}
\LABEL{s:Oval}

In this section
we investigate an oval of label $m$
for a minimal chart $\Gamma$.

Let $\Gamma$ be a chart, $m$ a label of $\Gamma$, $D$ a disk with $\partial D\subset \Gamma_m$, 
and $k$ a positive integer.
If $\partial D$ contains exactly
$k$ white vertices, 
then $D$ is called 
{\it a $k$-angled disk of $\Gamma_m$}. 
Note that 
the boundary $\partial D$ may contain crossings.

\begin{lemma}
\LABEL{LemmaOvalLens}
Let $\Gamma$ be a minimal chart, and
$m$ a label of $\Gamma$.
Let $G$ be an oval of label $m$, and
 $D$ a $2$-angled disk of $\Gamma_m$
with $\partial D\subset G$.
Let $E$ be a disk in $D$ whose boundary consists of
an internal edge in $G$ and
 an internal edge of label $m+\varepsilon$
$(\varepsilon\in\{+1,-1\})$
connecting the two white vertices
of $G$.
If $E$ does not contain
the terminal edges in $G$,
then 
$E$ is a lens of $\Gamma$.
\end{lemma}

\begin{Proof}
Let $e$ be the internal edge of label $m+\varepsilon$
in $\partial E$.
Let $v_1,v_2$ be the white vertices in $G$, and
$e_1,e_2$ the terminal edges at $v_1,v_2$ in $G$
respectively.
By Assumption~\ref{AssumeTerminal},
\begin{enumerate}
\item[$(1)$] both of the two edges $e_1$ and $e_2$
contain middle arcs.
\end{enumerate}
There are three cases:
(i) neither $e_1$ nor $e_2$ is contained in $D$
(see Fig.~\ref{Fig11}(a)),
(ii) only one of $e_1$ and $e_2$ is contained in $D$
(see Fig.~\ref{Fig11}(b)),
(iii) both of $e_1$ and $e_2$ are contained in $D$
(see Fig.~\ref{Fig11}(c)).

{\bf Case (i).}
By (1),
the edge $e$ is middle at 
both white vertices $v_1$ and $v_2$
simultaneously.
Thus the disk $E$ is a lens.

{\bf Case (ii).}
Without loss of generality
we can assume that 
\begin{enumerate}
\item[$(2)$] the edge $e_1$ is oriented inward at $v_1$.
\end{enumerate}
Then by (1),
the other two internal edges in $G$
are oriented from $v_1$ to $v_2$.
Thus
\begin{enumerate}
\item[$(3)$] the edge $e_2$ is oriented outward at $v_2$.
\end{enumerate}

If $e_1\subset D$ (see Fig.~\ref{Fig11}(b)),
then by (1) and (2),
the edge $e$ is oriented from $v_2$ to $v_1$.
Hence $e$ is oriented outward at $v_2$.
On the other hand,
by (1) and (3),
the edge $e$ is oriented inward at $v_2$.
This is a contradiction.

Similarly 
if $e_2\subset D$,
then we have the same contradiction.
Thus Case (ii) does not occur.

{\bf Case (iii).}
By (1),
none of $e$ and the two internal edges in $G$
contain middle arcs.
Thus the disk $E$ is a lens.
\end{Proof}

%%%%%%%%%%%%%%%%%%
%%%%%%%%%%%%%%%%%% Figure
%%%%%%%%%%%%%%%%%%
\begin{figure}[htb]
\centerline{\includegraphics{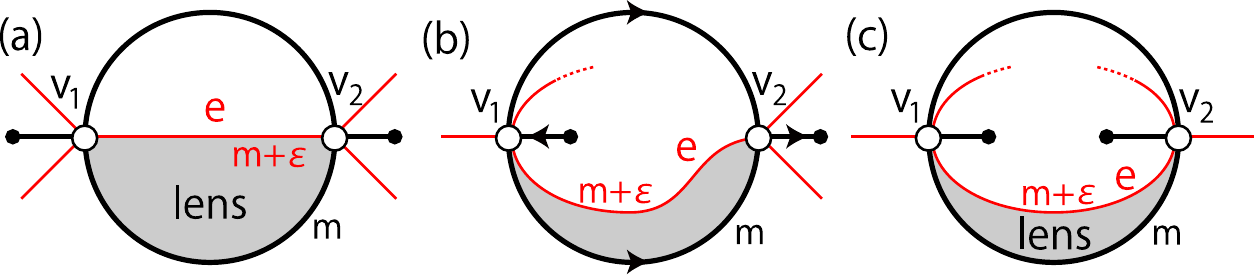}}
\caption{\LABEL{Fig11}
The gray regions are disks $E$.
(a) $e_1\not\subset D,e_2\not\subset D$.
(b) $e_1\subset D,e_2\not\subset D$.
(c) $e_1\subset D,e_2\subset D$.}
\end{figure}

Let $\Gamma$ be a chart. 
Suppose that an object consists of 
some edges of $\Gamma$, arcs in edges of $\Gamma$ and arcs around white vertices.
Then the object is called a {\it pseudo chart}.

\begin{lemma}
\LABEL{LemmaOval2angledDisk}
Let $\Gamma$ be a minimal chart, and
$m$ a label of $\Gamma$.
Let $G$ be an oval of label $m$.
If for some $\varepsilon\in\{+1,-1\}$
there exists a $2$-angled disk $D$ of 
$\Gamma_{m+\varepsilon}$ 
with $G\cap \partial D$
two white vertices,
then 
there exists a lens of $\Gamma$.
\end{lemma}

\begin{Proof}
Let $C$ be the simple closed curve in $G$.
Let $c_1,c_2,d_1,d_2$ be internal edges
with $c_1\cup c_2=C$ and $d_1\cup d_2=\partial D$.
Let $w$ be a white vertex in $G$.
Then there are two cases:
(i) The four edges $c_1,c_2,d_1,d_2$ (or $c_1,c_2,d_2,d_1$)
lie around the white vertex $w$ in this order
(see Fig~\ref{Fig12}(a)).
(ii) The four edges $c_1,d_1,c_2,d_2$
lie around the white vertex $w$ in this order
(see Fig~\ref{Fig12}(b)).

{\bf Case (i).}
Since $G\cap \partial D=C\cap \partial D$
consists of two white vertices,
the union $C\cup \partial D$
separates $S^2$ into four disks.
Let $E,E'$ be two of the four disks
such that each boundary consists of 
an internal edge of label $m$ in $G$ and
an internal edge of label $m+\varepsilon$ 
in $\partial D$.
By the condition of Case (i),
neither $E$ nor $E'$ contains a terminal edge
of label $m$ in $G$
(see Fig~\ref{Fig12}(a)).
Thus by Lemma~\ref{LemmaOvalLens},
both of $E$ and $E'$ are lenses.

{\bf Case (ii).}
Since $G\cap \partial D=C\cap \partial D$
consists of two white vertices,
the union $C\cup \partial D$
separates $S^2$ into four disks.
Let $E,E'$ be two of the four disks
with $E\cap E'=d_1$.
Then $E\cup E'$ is a 2-angled disk of $\Gamma_m$
whose boundary is contained in $G$.

If $E$ or $E'$ is a lens,
then there exists a lens of $\Gamma$.

Suppose that neither $E$ nor $E'$ is a lens.
Then by Lemma~\ref{LemmaOvalLens},
\begin{enumerate}
\item[$(1)$]
each of $E$ and $E'$ contains a terminal edge in $G$.
\end{enumerate}
Hence $E\cup E'$ contains one of 
the two pseudo charts as shown 
in Fig.~\ref{Fig12}(c),(d).
Then $d_2\not\subset E\cup E'$.
Thus $d_2\subset Cl(S^2-(E\cup E'))$.
Hence the edge $d_2$ separates 
the disk $Cl(S^2-(E\cup E'))$ into two disks.
Thus 
neither of the two disks 
contains a terminal edge in $G$,
because $E\cup E'$ contains the two terminal edges of $G$
by (1).
Hence by Lemma~\ref{LemmaOvalLens},
both of the two disks are lenses.
\end{Proof}

%%%%%%%%%%%%%%%%%%
%%%%%%%%%%%%%%%%%% Figure
%%%%%%%%%%%%%%%%%%
\begin{figure}[htb]
\centerline{\includegraphics{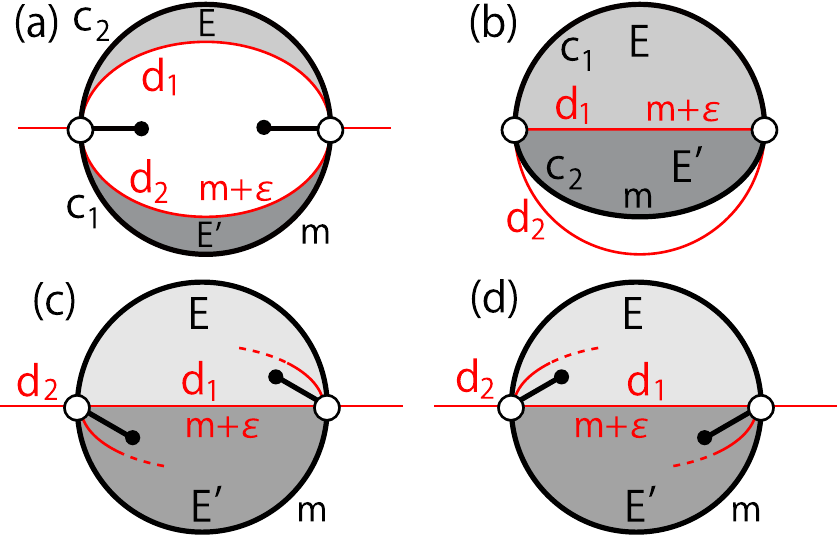}}
\caption{\LABEL{Fig12}
The light gray region and the dark gray region
are $E$ and $E'$.}
\end{figure}

\begin{lemma}
\LABEL{LemmaOvalNoGraph}
Let $\Gamma$ be a chart, and
$m$ a label of $\Gamma$.
Let $G$ be an oval of label $m$, and
$v_1,v_2$ the white vertices in $G$.
Let $D$ be a $2$-angled disk of $\Gamma_m$
with $\partial D\subset G$.
If $D$ satisfies one of the following two conditions,
then 
 $\Gamma$ is not minimal.
\begin{enumerate}
\item[{\rm (a)}] The disk $D$ 
does not contain terminal edges of $G$, but
contains two internal edges $e_1,e_2$ of 
label $m+\varepsilon$ at $v_1,v_2$ respectively
$(\varepsilon\in\{+1,-1\})$
such that
$e_1\cap e_2$ is a BW-vertex
with respect to $\Gamma_{m+\varepsilon}$ in ${\rm Int}D$
$($see Fig.~\ref{Fig13}$($a$))$.
\item[{\rm (b)}]
The disk $D$ contains exactly one terminal edge of $G$, 
and
contains three internal edges $e_1,e_2,e_3$ 
of label $m+\varepsilon$ at $v_1,v_1,v_2$
respectively
$(\varepsilon\in\{+1,-1\})$
such that 
$e_1\cap e_2\cap e_3$ is a white vertex
in ${\rm Int}D$
$($see Fig.~\ref{Fig13}$($b$))$.
\end{enumerate}
\end{lemma}

%%%%%%%%%%%%%%%%%%
%%%%%%%%%%%%%%%%%% Figure
%%%%%%%%%%%%%%%%%%
\begin{figure}[htb]
\centerline{\includegraphics{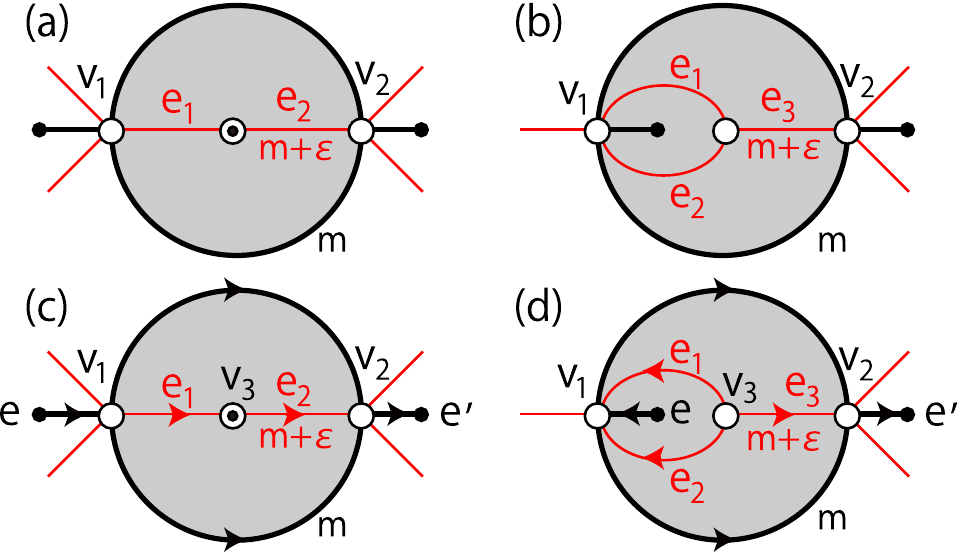}}
\caption{\LABEL{Fig13}
(a),(c) The gray regions are 2-angled disks 
not containing
terminal edges of $G$. 
(b),(d) The gray regions are 2-angled disks containing
one terminal edge of $G$.}
\end{figure}

\begin{Proof}
Suppose that
$\Gamma$ is minimal.
Let $e,e'$ be the terminal edges at $v_1,v_2$
in $G$ respectively.
By Assumption~\ref{AssumeTerminal},
\begin{enumerate}
\item[$(1)$] each of the two edges $e$ and $e'$
is middle at a white vertex.
\end{enumerate}
Without loss of generality
we can assume that 
\begin{enumerate}
\item[$(2)$] the edge $e$ is oriented inward at $v_1$.
\end{enumerate}
Then by (1),
the other two internal edges in $G$
are oriented from $v_1$ to $v_2$.
Thus
\begin{enumerate}
\item[$(3)$] the edge $e'$ is oriented outward at $v_2$.
\end{enumerate}

If the disk $D$ satisfies the condition (a),
then $e_1\cap e_2$ is a BW-vertex with respect to 
$\Gamma_{m+\varepsilon}$.
Let $v_3=e_1\cap e_2$.
By (2),
the edge $e_1$ is oriented from $v_1$ to $v_3$.
Thus 
the edge $e_1$ is oriented inward at the BW-vertex $v_3$.

On the other hand, by (3),
the edge $e_2$ is oriented from $v_3$ to $v_2$.
Thus 
the edge $e_2$ is oriented outward at the BW-vertex $v_3$.
Hence for the BW-vertex $v_3$,
the edge $e_1$ of label $m+\varepsilon$
is oriented inward at $v_3$,
and the edge $e_2$ of label $m+\varepsilon$
is oriented outward at $v_3$
(see Fig.~\ref{Fig13}(c)).
This contradicts Lemma~\ref{OriBWvertex}.
Thus $\Gamma$ is not minimal.

If the disk $D$ satisfies the condition (b),
then $e_1\cap e_2\cap e_3$ is a white vertex.
Let $v_3=e_1\cap e_2\cap e_3$.
By (1) and (2),
both of $e_1$ and $e_2$ are oriented from $v_3$ to $v_1$.
Thus 
\begin{enumerate}
\item[$(4)$] 
both of $e_1$ and $e_2$ are oriented outward at $v_3$.
\end{enumerate}

On the other hand, by (3),
the edge $e_3$ is oriented from $v_3$ to $v_2$.
Thus 
the edge $e_3$ is oriented outward at $v_3$.
Hence by (4),
the three edges $e_1,e_2,e_3$ of label $m+\varepsilon$
are oriented outward at $v_3$
(see Fig.~\ref{Fig13}(d)).
This contradicts the condition (iii) of 
the definition of charts.
Thus $\Gamma$ is not minimal.
\end{Proof}

%%%%%%%%%%%%%%%%%%%%%
%%%%%%%%%%%%%%%%%%%%%

%\newpage

\section{White vertices in the graphs $\Gamma_{m+1}$}
\LABEL{s:WhiteGammaM+1GammaM+2}

In this section,
we investigate white vertices 
in an oval of label $m$ 
for a minimal chart $\Gamma$ of type $(m;2,3,2)$.

\begin{lemma}
\LABEL{LemmaTwoEdges}
Let $\Gamma$ be a minimal chart,
and $m$ a label of $\Gamma$.
Let $w$ be a white vertex in 
a terminal edge of label $m$.
Let $e_1,e_2$ be the two edges of label $m$ at $w$
different from the terminal edge.
If $w$ is contained in a terminal edge of label
$m+\varepsilon$ for some $\varepsilon\in\{+1,-1\}$,
then 
both edges $e_1,e_2$ are contained in
the closure of the same connected component 
of $S^2-\Gamma_{m+\varepsilon}$.
\end{lemma}

\begin{Proof}
Since $w$ is contained in a terminal edge of label $m$ and 
since $w$ is contained in a terminal edge of label $m+\varepsilon$,
by Assumption~\ref{AssumeTerminal}
we can show
that 
in a neighborhood of the vertex $w$,
the chart $\Gamma$ contains the pseudo chart 
as shown in Fig.~\ref{Fig14}.
Hence 
the edges $e_1$ and $e_2$
of label $m$ 
are contained in
the closure of the same 
connected component $F$ of $S^2-\Gamma_{m+\varepsilon}$.
Thus we complete the proof of 
Lemma~\ref{LemmaTwoEdges}.
\end{Proof}

%%%%%%%%%%%%%%%%%%
%%%%%%%%%%%%%%%%%% Figure
%%%%%%%%%%%%%%%%%%
\begin{figure}[htb]
\centerline{\includegraphics{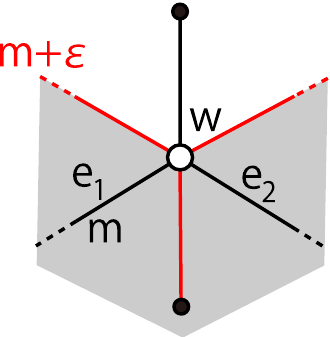}}
\caption{\LABEL{Fig14}
A white vertex $w$ is contained in two terminal edges.
The gray region is $F$.}
\end{figure}

From the above lemma,
we have the following lemma:

\begin{lemma}
\LABEL{LemmaOvalTwoEdges}
Let $\Gamma$ be a minimal chart,
and $m$ a label of $\Gamma$.
Let $G$ be an oval of label $m$.
If one of the two white vertices in $G$ is a BW-vertex with respect to $\Gamma_{m+\varepsilon}$ for some 
$\varepsilon\in\{+1,-1\}$,
then the two internal edges in $G$
are contained in
the closure of the same connected component 
of $S^2-\Gamma_{m+\varepsilon}$.
\end{lemma}

\begin{lemma}
\LABEL{LemmaBWvertexOval}
Let $\Gamma$ be a minimal chart of type $(m;2,3,2)$.
Suppose that $\Gamma_m$ contains an oval $G$.
Then we have the following.
\begin{enumerate}
\item[{\rm (a)}] 
If $\Gamma_{m+1}$ 
contains one of the three graphs as shown in 
Fig.~\ref{Fig04}$($a$)$,$($b$)$,$($c$)$,
then 
either $G$ contains two BW-vertices
with respect to $\Gamma_{m+1}$,
or  
$G$ does not contain any BW-vertex
with respect to $\Gamma_{m+1}$.
\item[{\rm (b)}]
If $\Gamma_{m+1}$ 
contains one of the three graphs as shown in 
Fig.~\ref{Fig04}$($d$)$,$($e$)$,$($f$)$,
then $G$ does not contain any BW-vertex 
with respect to $\Gamma_{m+1}$.
\item[{\rm (c)}]
If $\Gamma_{m+1}$
contains the union of a $\theta$-curve and 
a skew $\theta$-curve,
then 
 $G$ does not contain any BW-vertex 
with respect to $\Gamma_{m+1}$.
Moreover 
the two white vertices in $G$ 
are contained in the $\theta$-curve
or the skew $\theta$-curve simultaneously.
\item[{\rm (d)}]
If $\Gamma_{m+1}$ contains the union of an oval
and a skew $\theta$-curve,
then either $G$ contains two BW-vertices
with respect to $\Gamma_{m+1}$,
or  
$G$ does not contain any BW-vertex
with respect to $\Gamma_{m+1}$.
\end{enumerate}
\end{lemma}

\begin{Proof}
Let $e_1,e_2$ be the internal edges of label $m$
in $G$.

{\bf Statement (a).}
The graph $\Gamma_{m+1}$
contains exactly three  BW-vertices
with respect to $\Gamma_{m+1}$, 
say $w_1,w_2,w_3$.
Let $w_4,w_5$ be
the other white vertices in $\Gamma_{m+1}$.
It suffices to prove that
if $G$ contains one of BW-vertices $w_1,w_2,w_3$,
then $G$ contains two of $w_1,w_2,w_3$.

If $G$ contains one of BW-vertices $w_1,w_2,w_3$,
then by Lemma~\ref{LemmaOvalTwoEdges}
there exists a connected component $F$ of 
$S^2-\Gamma_{m+1}$
with $e_1\cup e_2\subset Cl(F)$
(see Fig.~\ref{Fig15}(a)).
Thus 
\begin{enumerate}
\item[(1)] 
for each white vertex $w$ in $G$,
there exist {\it two} edges of label $m$ 
at $w$ contained in $Cl(F)$.
\end{enumerate}

On the other hand,
by the condition of Lemma~\ref{LemmaBWvertexOval}(a),
for each white vertex $w_i$ ($i=4,5$)
there exists {\it at most one} edge label $m$ at $w_i$
in $Cl(F)$ (see Fig.~\ref{Fig15}(a)).
Hence by (1),
the oval $G$ does not contain 
$w_4$ nor $w_5$.
Thus $G$ contains two of BW-vertices $w_1,w_2,w_3$.
Hence Statement (a) holds.

{\bf Statement (b).}
The graph $\Gamma_{m+1}$ contains 
exactly one BW-vertex
with respect to $\Gamma_{m+1}$, 
say $w_1$.
Let $w_2,w_3,w_4,w_5$ be
the other white vertices in $\Gamma_{m+1}$.

Suppose that $G$ contains the BW-vertex $w_1$.
Then 
by Lemma~\ref{LemmaOvalTwoEdges}
there exists a connected component $F$ of $S^2-\Gamma_{m+1}$ with $e_1\cup e_2\subset Cl(F)$
(see Fig.~\ref{Fig15}(b)).
Thus 
for each white vertex in the oval $G$,
there exist {\it two} edges of label $m$ 
at the vertex in $Cl(F)$.
However, 
by the condition of Lemma~\ref{LemmaBWvertexOval}(b),
for each white vertex $w_i$
$(i=2,3,4,5)$
there exists {\it at most one} edge of label $m$ at $w_i$
in $Cl(F)$.
This is a contradiction.
Hence the oval $G$ does not contain 
the vertex $w_1$.
Hence Statement (b) holds.

{\bf Statement (c).}
Let $w_1,w_2$ be the white vertices 
of the $\theta$-curve in $\Gamma_{m+1}$.
Let $w_3$ be the BW-vertex 
of the skew $\theta$-curve
with respect to $\Gamma_{m+1}$,
and $w_4,w_5$ the other white vertices 
of the skew $\theta$-curve.

By a similar way of the proof of Statement (b),
we can show that
the oval $G$ does not contain the BW-vertex $w_3$.
Thus $G$ contains two of $w_1,w_2,w_4,w_5$.

Suppose that the oval $G$ contains 
one of the white vertices $w_1,w_2$,
and one of the white vertices $w_4,w_5$.
Without loss of generality
we can assume $w_1,w_4\in G$.
Then 
\begin{enumerate}
\item[$(2)$]
the two edges $e_1$ and $e_2$ of label $m$
connect the vertices $w_1$ and $w_4$.
\end{enumerate}

Now, the $\theta$-curve in $\Gamma_{m+1}$
separates $S^2$ into three disks.
One contains the skew $\theta$-curve in $\Gamma_{m+1}$,
say $F$.
Hence the vertex $w_4$ in the skew $\theta$-curve
is contained in $F$.
Thus by (2),
we have $e_1\cup e_2\subset F$.
Hence 
in $F$
there exist {\it two} edges of label $m$ at $w_1$.
However,
since $w_1$ is a white vertex 
of the $\theta$-curve in $\Gamma_{m+1}$,
there exists {\it at most one} edge of label $m$ at $w_1$
in $F$ (see Fig.~\ref{Fig15}(c)).
This is a contradiction.
Hence $w_1,w_2\in G$ or
$w_4,w_5\in G$.
Thus Statement (c) holds.

Similarly we can show  Statement (d).
\end{Proof}

%%%%%%%%%%%%%%%%%%
%%%%%%%%%%%%%%%%%% Figure
%%%%%%%%%%%%%%%%%%
\begin{figure}[htb]
\centerline{\includegraphics{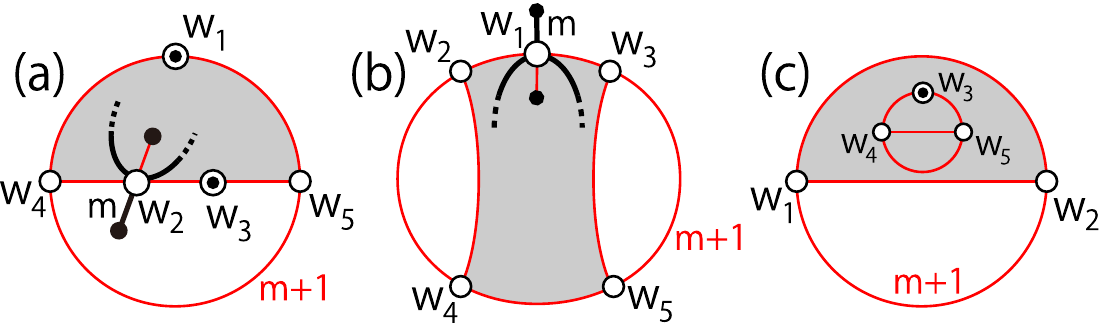}}
\caption{\LABEL{Fig15}
The gray regions are $F$.
(a) The graph as shown in Fig.~\ref{Fig04}(b)
with $w_2\in G$.
(b) The graph as shown in Fig.~\ref{Fig04}(d)
with $w_1\in G$.
(c) The skew $\theta$-curve of label $m+1$ 
is contained in $F$.
}
\end{figure}

%%%%%%%%%%%%%%%%%%%%%%%%%%%%%%%%%%%%%%%%%
%%%%%%%%%%%%%%%%%%%%%%%%%%%%%%%%%%%%%%%%%
%%%%%%%%%%%%%%%%%%%%%%%%%%%%%%%%%%%%%%%%
%\newpage

\section{The graphs $\Gamma_{m+1}$ and $\Gamma_{m+2}$}
\LABEL{s:GammaM+1GammaM+2}

In this section,
we shall show that
for any minimal chart $\Gamma$ of type $(m;2,3,2)$,
the graph $\Gamma_{m+1}$ contains none of
the five graphs as shown in 
Fig.~\ref{Fig04}(a),(d),(e),(f),(i),
and neither does $\Gamma_{m+2}$.
Moreover
we shall show that
neither $\Gamma_{m+1}$ nor $\Gamma_{m+2}$
contains a $\theta$-curve.

\begin{lemma}
\LABEL{GammaM+1ToGammaM+2}
Let $G$ be one of 12 graphs as shown in
Fig.~\ref{Fig03} and Fig.~\ref{Fig04}.
If 
for any minimal chart $\Gamma$ of type $(m;2,3,2)$,
the graph $\Gamma_{m+1}$ does not contain
the graph $G$,
then the graph $\Gamma_{m+2}$ does not contain
the graph $G$.
\end{lemma}

\begin{Proof}
Suppose that the graph $\Gamma_{m+1}$
does not contain the graph $G$ 
for any minimal chart $\Gamma$ of type $(m;2,3,2)$.

If there exists a minimal chart $\Gamma'$ 
of type $(m;2,3,2)$ with $\Gamma_{m+2}'\supset G$,
then let $\Gamma''$ be the chart obtained from $\Gamma'$
by changing labels $\cdots,m,m+1,m+2,m+3,\cdots$
into $\cdots,m+3,m+2,m+1,m,\cdots$,
respectively.
Then $\Gamma''$ is a chart of type $(m;2,3,2)$
with $\Gamma_{m+1}''\supset G$.
Hence $\Gamma''$ is not minimal.
Thus $\Gamma''$ is C-move equivalent to 
a chart whose complexity is less than 
the complexity of $\Gamma''$.
Hence by using the above C-moves,
the chart $\Gamma'$ is also C-move equivalent to 
a chart whose complexity is less than 
the complexity $\Gamma'$.
Thus $\Gamma'$ is not minimal.
This is a contradiction.
Hence 
if $\Gamma'$ is a minimal chart of type $(m;2,3,2)$,
then $\Gamma_{m+2}'\not\supset G$.
\end{Proof}

\begin{lemma}
\LABEL{LemmaTypeA}
Let $\Gamma$ be a minimal chart of type $(m;2,3,2)$.
Then neither $\Gamma_{m+1}$ nor $\Gamma_{m+2}$
contains the graph as shown in 
Fig.~\ref{Fig04}$($a$)$.
\end{lemma}

\begin{Proof}
Suppose that $\Gamma_{m+1}$
contains the graph as shown
in Fig.~\ref{Fig04}$($a$)$.
We use the notations as shown in Fig.~\ref{Fig16}(a)
where $w_3,w_4,w_5$ are BW-vertices.
By Corollary~\ref{NoThetaGammaM},
the graph $\Gamma_m$ contains an oval $G$.
Thus by Lemma~\ref{LemmaBWvertexOval}(a),
there are three cases:
(i) $w_1,w_2\in G$,
(ii) $w_3,w_4\in G$ or $w_4,w_5\in G$
(see Fig.~\ref{Fig16}(b)),
(iii) $w_3,w_5\in G$
(see Fig.~\ref{Fig16}(c)).

{\bf Case (i).}
Since there exist two internal edges of 
label $m+1$ connecting $w_1$ and $w_2$,
there exists a 2-angled disk $D$ of $\Gamma_{m+1}$
with $w_1,w_2\in\partial D$.
Thus $w_1,w_2\in G\cap \partial D$.
Hence by Lemma~\ref{LemmaOval2angledDisk},
there exists a lens of $\Gamma$.
This contradicts Lemma~\ref{NoLens}.
Hence Case (i) does not occur.

{\bf Case (ii).}
By Lemma~\ref{LemmaOvalTwoEdges},
the two internal edges $e_1,e_2$ of label $m$ in $G$
are contained in the closure of 
the same connected component $F$ of $S^2-\Gamma_{m+1}$.
Thus the curve $e_1\cup e_2$ bounds 
a 2-angled disk of $\Gamma_m$ in $Cl(F)$,
say $D$.
Hence $Cl(S^2-D)$ is also a 2-angled disk of $\Gamma_m$,
and by Lemma~\ref{LemmaOvalLens}
the disk $Cl(S^2-D)$ contains a lens 
(see Fig.~\ref{Fig16}(b)).
This contradicts Lemma~\ref{NoLens}.
Thus Case (ii) does not occur.

{\bf Case (iii).}
Let $e_1,e_2$ be the two internal edges 
of label $m$ in $G$.
Let $e_3,e_5$ be the internal edges of label $m+1$
at $w_3,w_5$ containing $w_4$ respectively
(see Fig.~\ref{Fig16}(c)).

By Lemma~\ref{LemmaOvalTwoEdges},
the two edges $e_1,e_2$ are contained in the closure of 
the same connected component $F$ of $S^2-\Gamma_{m+1}$
(see Fig.~\ref{Fig16}(c)).
 Without loss of generality
we can assume that
the terminal edge of label $m$ at $w_3$ 
is oriented inward at $w_3$.
By Assumption~\ref{AssumeTerminal},
the terminal edge is middle at $w_3$.
Thus 
\begin{enumerate}
\item[$(1)$]
the edge $e_3$ is oriented inward at $w_3$,
\end{enumerate}
and the two edges $e_1,e_2$ 
are oriented from $w_3$ to $w_5$.
Hence the edge $e_5$ is
oriented outward at $w_5$
(see Fig~\ref{Fig16}(d)).
Thus $e_5$
is oriented inward at the BW-vertex $w_4$.
However by (1)
the edge $e_3$
is oriented outward at the BW-vertex $w_4$.
This contradicts Lemma~\ref{OriBWvertex}.
Hence Case (iii) does not occur.

Therefore
 $\Gamma_{m+1}$ does not
contain the graph as shown in 
Fig.~\ref{Fig04}$($a$)$.
By Lemma~\ref{GammaM+1ToGammaM+2}, 
we can show that
 $\Gamma_{m+2}$ does not
contain the graph as shown in 
Fig.~\ref{Fig04}$($a$)$.
\end{Proof}

%%%%%%%%%%%%%%%%%%
%%%%%%%%%%%%%%%%%% Figure
%%%%%%%%%%%%%%%%%%
\begin{figure}[htb]
\centerline{\includegraphics{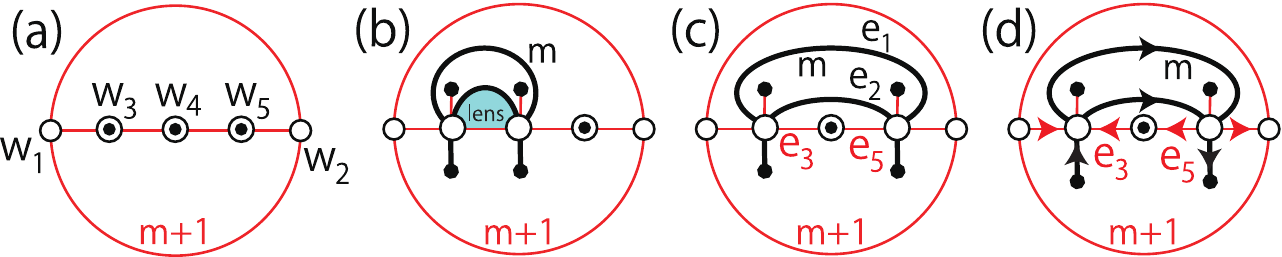}}
\caption{\LABEL{Fig16}
(a) $w_1,w_2,\cdots,w_5$ are white vertices.
(b) $w_3,w_4\in G$.
(c),(d) $w_3,w_5\in G$.}
\end{figure}

%%%%%%%%%%%%%%%%%%%%%%%%%%%%

\begin{lemma}
\LABEL{LemmaTypeD}
Let $\Gamma$ be a minimal chart of type $(m;2,3,2)$.
Then neither $\Gamma_{m+1}$ nor $\Gamma_{m+2}$
contains the graph as shown in 
Fig.~\ref{Fig04}$($d$)$.
\end{lemma}

\begin{Proof}
Suppose that $\Gamma_{m+1}$
contains the graph as shown 
in Fig.~\ref{Fig04}$($d$)$.
We use the notations as shown in Fig.~\ref{Fig17}(a)
where $w_1$ is a BW-vertex.
By Corollary~\ref{NoThetaGammaM},
the graph $\Gamma_m$ contains an oval $G$.
Thus by Lemma~\ref{LemmaBWvertexOval}(b),
 there are four cases:
(i) $w_2,w_3\in G$ 
(see Fig.~\ref{Fig17}(b)),
(ii) $w_2,w_4\in G$ or $w_3,w_5\in G$,
(iii) $w_2,w_5\in G$ or $w_3,w_4\in G$
(see Fig.~\ref{Fig17}(c)),
(iv) $w_4,w_5\in G$
(see Fig.~\ref{Fig17}(d)).

{\bf Case (i).}
By Lemma~\ref{LemmaOvalNoGraph}(a),
the chart $\Gamma$ is not minimal.
This is a contradiction.
Hence Case (i) does not occur.

{\bf Case (ii).}
By Lemma~\ref{LemmaOval2angledDisk},
there exists a lens of $\Gamma$.
This contradicts Lemma~\ref{NoLens}.
Hence Case (ii) does not occur.

{\bf Case (iii).}
By Lemma~\ref{LemmaOvalNoGraph}(b),
the chart $\Gamma$ is not minimal.
This is a contradiction.
Hence Case (iii) does not occur.

{\bf Case (iv).}
By Lemma~\ref{LemmaOvalLens},
there exists a lens of $\Gamma$.
This contradicts Lemma~\ref{NoLens}.
Hence Case (iv) does not occur.

Therefore
 $\Gamma_{m+1}$ does not
contain the graph as shown in 
Fig.~\ref{Fig04}$($d$)$.
By Lemma~\ref{GammaM+1ToGammaM+2}, 
we can show that
 $\Gamma_{m+2}$ does not
contain the graph as shown in 
Fig.~\ref{Fig04}$($d$)$.
\end{Proof}

%%%%%%%%%%%%%%%%%%
%%%%%%%%%%%%%%%%%% Figure
%%%%%%%%%%%%%%%%%%
\begin{figure}[htb]
\centerline{\includegraphics{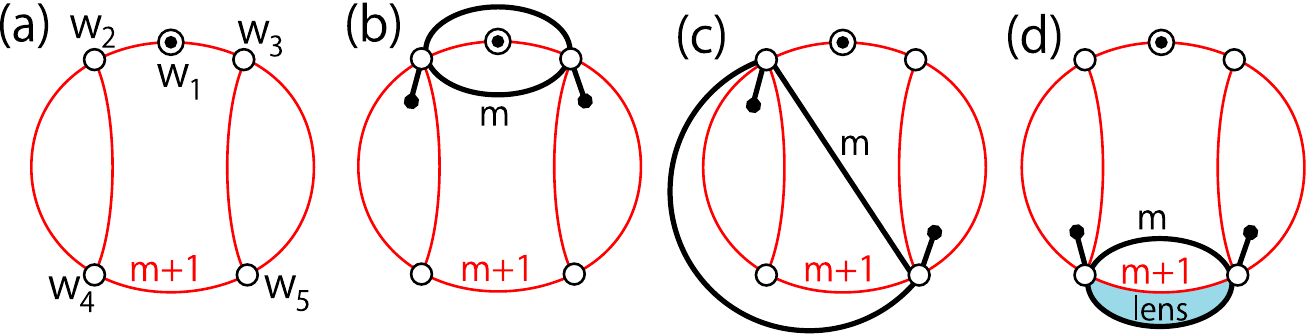}}
\caption{\LABEL{Fig17}
(a) $w_1,w_2,\cdots,w_5$ are white vertices.
(b) $w_2,w_3\in G$.
(c) $w_2,w_5\in G$.
(d)  $w_4,w_5\in G$.}
\end{figure}

%%%%%%%%%%%%%%%%%%%%%%%%%%%%
%%%%%%%%%%%%%%%%%%%%%%%%%%%%
%%%%%%%%%%%%%%%%%%%%%%%%%%%%
%%%%%%%%%%%%%%%%%%%%%%%%%%%%

\begin{lemma}
\LABEL{LemmaTypeE}
Let $\Gamma$ be a minimal chart of type $(m;2,3,2)$.
Then neither $\Gamma_{m+1}$ nor $\Gamma_{m+2}$
contains the graph as shown in 
Fig.~\ref{Fig04}$($e$)$.
\end{lemma}

\begin{Proof}
Suppose that $\Gamma_{m+1}$
contains the graph as shown in 
Fig.~\ref{Fig04}$($e$)$.
We use the notations as shown in Fig.~\ref{Fig18}(a)
where $w_1$ is a BW-vertex.
By Corollary~\ref{NoThetaGammaM},
the graph $\Gamma_m$ contains an oval $G$.
Thus by Lemma~\ref{LemmaBWvertexOval}(b),
there are four cases:
(i) $w_2,w_3\in G$ 
(see Fig.~\ref{Fig18}(b),(c),(d)),
(ii) $w_2,w_4\in G$ or $w_3,w_5\in G$
(see Fig.~\ref{Fig18}(e)),
(iii) $w_2,w_5\in G$ or $w_3,w_4\in G$
(see Fig.~\ref{Fig18}(f)),
(iv) $w_4,w_5\in G$.

{\bf Case (i).}
Around the white vertex $w_2$,
there are three internal edges of label $m+1$.
Let $e_1,e_2,e_3$ be the three internal edges
of label $m+1$ at $w_2$
containing $w_1,w_3,w_4$,
respectively
(see Fig.~\ref{Fig18}(a)).
Let $D$ be the 2-angled disk of $\Gamma_m$
not containing the terminal edge of label $m$ at $w_2$.
Then there are three cases:
$e_1\subset D$ or $e_2\subset D$
or $e_3\subset D$.

If $e_1\subset D$ (see Fig.~\ref{Fig18}(b)),
then by Lemma~\ref{LemmaOvalNoGraph}(a)
the chart $\Gamma$ is not minimal.
This is a contradiction.
If $e_2\subset D$  (see Fig.~\ref{Fig18}(c)),
then by Lemma~\ref{LemmaOvalLens}
there exists a lens of $\Gamma$.
This contradicts Lemma~\ref{NoLens}.
If $e_3\subset D$ (see Fig.~\ref{Fig18}(d)),
then by Lemma~\ref{LemmaOvalLens}
there exists a lens in $Cl(S^2-D)$.
This contradicts Lemma~\ref{NoLens}.
Hence Case (i) does not occur.

{\bf Case (ii).}
By Lemma~\ref{LemmaOvalLens},
there exists a lens of $\Gamma$.
This contradicts Lemma~\ref{NoLens}.
Hence Case (ii) does not occur.

{\bf Case (iii).}
By Lemma~\ref{LemmaOvalNoGraph}(b),
the chart $\Gamma$ is not minimal.
This is a contradiction.
Hence Case (iii) does not occur.

{\bf Case (iv).}
By Lemma~\ref{LemmaOval2angledDisk},
there exists a lens of $\Gamma$.
This contradicts Lemma~\ref{NoLens}.
Hence Case (iv) does not occur.

Therefore
 $\Gamma_{m+1}$ does not
contain the graph as shown in 
Fig.~\ref{Fig04}$($e$)$.
By Lemma~\ref{GammaM+1ToGammaM+2}, 
we can show that
 $\Gamma_{m+2}$ does not
contain the graph as shown in 
Fig.~\ref{Fig04}$($e$)$.
\end{Proof}

%%%%%%%%%%%%%%%%%%
%%%%%%%%%%%%%%%%%% Figure
%%%%%%%%%%%%%%%%%%
\begin{figure}[htb]
\centerline{\includegraphics{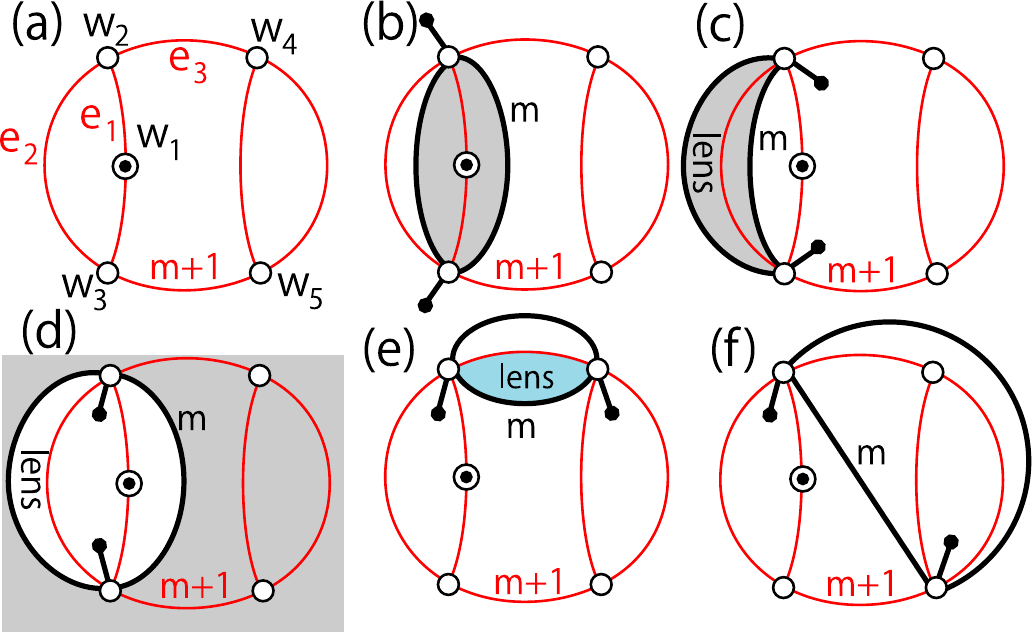}}
\caption{\LABEL{Fig18}
(a) $w_1,w_2,\cdots,w_5$ are white vertices.
(b), (c), (d) $w_2,w_3\in G$,
the gray regions are 2-angled disks of $\Gamma_m$
not containing the terminal edges of $G$.
(e) $w_2,w_4\in G$.
(f) $w_2,w_5\in G$.}
\end{figure}

\begin{lemma}
\LABEL{LemmaTypeF}
Let $\Gamma$ be a minimal chart of type $(m;2,3,2)$.
Then neither $\Gamma_{m+1}$ nor $\Gamma_{m+2}$
contains the graph as shown in 
Fig.~\ref{Fig04}$($f$)$.
\end{lemma}

\begin{Proof}
Suppose that $\Gamma_{m+1}$
contains the graph as shown in 
Fig.~\ref{Fig04}$($f$)$.
We use the notations as shown in Fig.~\ref{Fig19}(a)
where $w_1$ is a BW-vertex.
By Corollary~\ref{NoThetaGammaM},
the graph $\Gamma_m$ contains an oval $G$.
There are two cases:
(i) $w_4\in G$ or $w_5\in G$,
(ii) $w_4\not\in G$ and $w_5\not\in G$.

{\bf Case (i).}
If $w_4\in G$,
then by Lemma~\ref{LemmaBWvertexOval}(b)
the oval $G$ contains one of $w_2,w_3,w_5$.
Thus there exists an internal edge of label $m+1$
connecting the two white vertices of $G$
(see Fig.~\ref{Fig19}(b)).
Hence by Lemma~\ref{LemmaOvalLens},
there exists a lens of $\Gamma$.
This contradicts Lemma~\ref{NoLens}.

If $w_5\in G$,
then we have the same contradiction.
Hence Case (i) does not occur.

{\bf Case (ii).}
By Lemma~\ref{LemmaBWvertexOval}(b),
we have $w_2,w_3\in G$
(see Fig.~\ref{Fig19}(c)).
By Lemma~\ref{LemmaOvalNoGraph}(a),
the chart $\Gamma$ is not minimal.
This is a contradiction.
Hence Case (ii) does not occur.

Therefore
 $\Gamma_{m+1}$ does not
contain the graph as shown in 
Fig.~\ref{Fig04}$($f$)$.
By Lemma~\ref{GammaM+1ToGammaM+2}, 
we can show that
 $\Gamma_{m+2}$ does not
contain the graph as shown in 
Fig.~\ref{Fig04}$($f$)$.
\end{Proof}

%%%%%%%%%%%%%%%%%%
%%%%%%%%%%%%%%%%%% Figure
%%%%%%%%%%%%%%%%%%
\begin{figure}[htb]
\centerline{\includegraphics{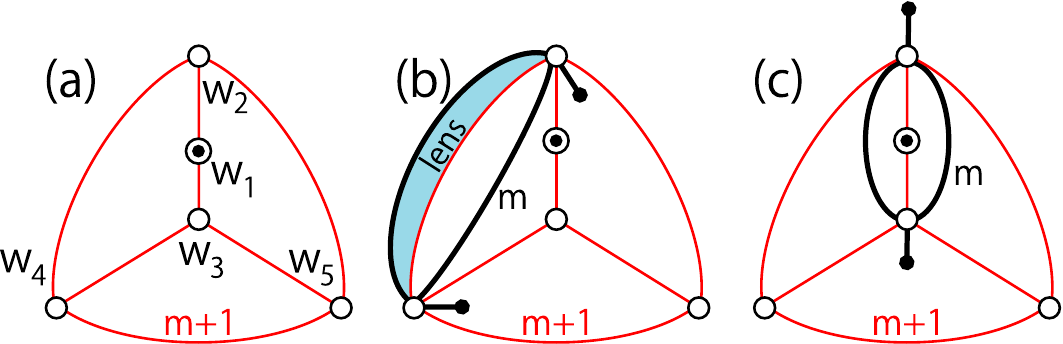}}
\caption{\LABEL{Fig19}
(a) $w_1,w_2,\cdots,w_5$ are white vertices.
(b) $w_2,w_4\in G$.
(c) $w_2,w_3\in G$.}
\end{figure}

\begin{lemma}
\LABEL{LemmaTypeI}
Let $\Gamma$ be a minimal chart of type $(m;2,3,2)$.
Then neither $\Gamma_{m+1}$ nor $\Gamma_{m+2}$
contains the graph as shown in 
Fig.~\ref{Fig04}$($i$)$.
\end{lemma}

\begin{Proof}
Suppose that $\Gamma_{m+1}$
contains the graph as shown in 
Fig.~\ref{Fig04}$($i$)$.
We use the notations as shown in Fig.~\ref{Fig20}(a)
where $w_1$ is a BW-vertex.
By Corollary~\ref{NoThetaGammaM},
the graph $\Gamma_m$ contains an oval $G$.
There are four cases:
(i) $G$ contains the BW-vertex $w_1$
 (see Fig.~\ref{Fig20}(b),(c)),
(ii) $w_2,w_3\in G$ 
(see Fig.~\ref{Fig20}(d)),
(iii) $w_3,w_4\in G$ or $w_3,w_5\in G$
(see Fig.~\ref{Fig20}(e)),
(iv) $w_4,w_5\in G$.

{\bf Case (i).}
Let $e_1,e_2$ be the internal edges in $G$,
and $F_1,F_2$ the connected components of
$S^2-\Gamma_{m+1}$
with $w_1\in Cl(F_1)\cap Cl(F_2)$
and $w_3\not\in Cl(F_1)$
(see Fig.~\ref{Fig20}(a)).

We shall show that $e_1\cup e_2\subset Cl(F_2)$.
By Lemma~\ref{LemmaOvalTwoEdges},
we have 
$e_1\cup e_2\subset Cl(F_1)$ or
$e_1\cup e_2\subset Cl(F_2)$.
If $e_1\cup e_2\subset Cl(F_1)$,
then in $Cl(F_1)$
there exist {\it two} edges of label $m$
at $w_2$.
However, 
since $w_2$ is a white vertex as shown 
in Fig.~\ref{Fig20}(a),
there exists {\it at most one} edge of label $m$
at $w_2$ in $Cl(F_1)$.
This is a contradiction.
Hence $e_1\cup e_2\subset Cl(F_2)$.

Thus there are two cases:
$w_1,w_2\in G$ or $w_1,w_3\in G$.
If $w_1,w_2\in G$ (see Fig.~\ref{Fig20}(b)),
then by Lemma~\ref{LemmaOvalLens}
there exists a lens of $\Gamma$.
This contradicts Lemma~\ref{NoLens}.
If $w_1,w_3\in G$ (see Fig.~\ref{Fig20}(c)),
then by Lemma~\ref{LemmaOvalNoGraph}(b)
the chart $\Gamma$ is not minimal.
This is a contradiction.
Hence Case (i) does not occur.

{\bf Case (ii) and Case (iii).}
By Lemma~\ref{LemmaOvalLens},
there exists a lens of $\Gamma$.
This contradicts Lemma~\ref{NoLens}.
Hence Case (ii) and Case (iii) do not occur.

{\bf Case (iv).}
By Lemma~\ref{LemmaOval2angledDisk},
there exists a lens of $\Gamma$.
This contradicts Lemma~\ref{NoLens}.
Hence Case (iv) does not occur.

Therefore
 $\Gamma_{m+1}$ does not
contain the graph as shown in 
Fig.~\ref{Fig04}$($i$)$.
By Lemma~\ref{GammaM+1ToGammaM+2}, 
we can show that
 $\Gamma_{m+2}$ does not
contain the graph as shown in 
Fig.~\ref{Fig04}$($i$)$.
\end{Proof}

%%%%%%%%%%%%%%%%%%
%%%%%%%%%%%%%%%%%% Figure
%%%%%%%%%%%%%%%%%%
\begin{figure}[htb]
\centerline{\includegraphics{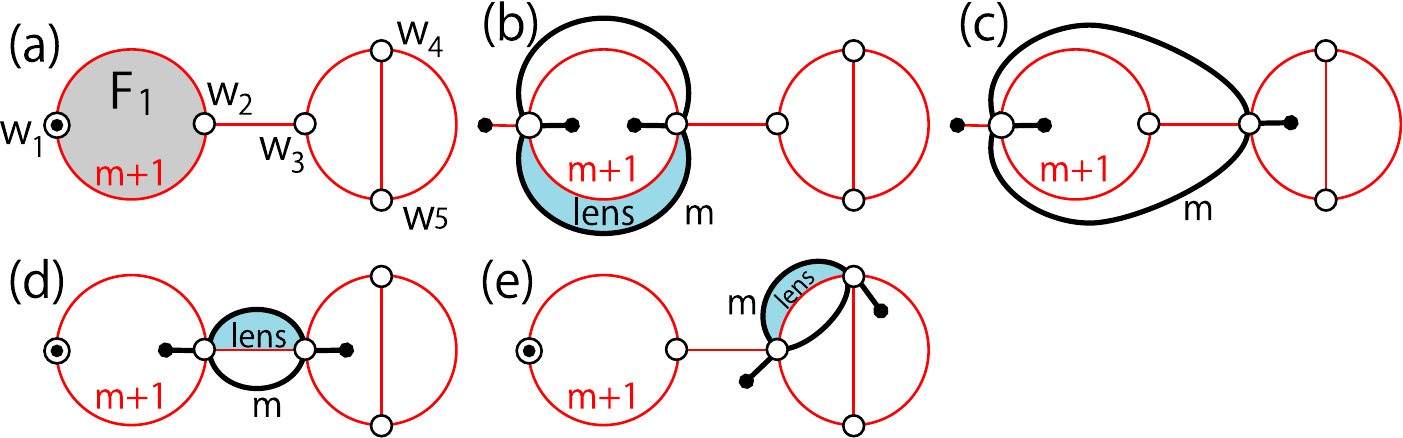}}
\caption{\LABEL{Fig20}
(a) $w_1,w_2,\cdots,w_5$ are white vertices,
the gray region is $F_1$.
(b) $w_1,w_2\in G$.
(c) $w_1,w_3\in G$.
(d) $w_2,w_3\in G$.
(e) $w_3,w_4\in G$.
}
\end{figure}

\begin{lemma}
\LABEL{TypeThetaSkewTheta}
Let $\Gamma$ be a minimal chart of type $(m;2,3,2)$.
Then neither $\Gamma_{m+1}$ nor $\Gamma_{m+2}$
contains a $\theta$-curve.
\end{lemma}

\begin{Proof}
Suppose that $\Gamma_{m+1}$
contains a $\theta$-curve $G_1$. 
By Lemma~\ref{GammaM+1},
the graph $\Gamma_{m+1}$
contains a skew $\theta$-curve $G_2$.
Let $w_1,w_2$ be the white vertices in $G_1$.
Let $w_3$ be the BW-vertex in $G_2$,
and $w_4,w_5$ the other white vertices in $G_2$.

By Corollary~\ref{NoThetaGammaM},
the graph $\Gamma_m$ contains an oval $G$.
Thus by Lemma~\ref{LemmaBWvertexOval}(c),
there are two cases:
 $w_1,w_2\in G$, or $w_4,w_5\in G$.
For both cases,
there exist two internal edges of label $m+1$
connecting the two white vertices in $G$.
Hence there exists a 2-angled disk $D$ of $\Gamma_{m+1}$
with $G\cap \partial D$ two white vertices.
Thus by Lemma~\ref{LemmaOval2angledDisk},
there exists a lens of $\Gamma$.
This contradicts Lemma~\ref{NoLens}.
Therefore
 $\Gamma_{m+1}$ does not contain a $\theta$-curve.

By Lemma~\ref{GammaM+1ToGammaM+2}, 
we can show that
 $\Gamma_{m+2}$ does not contain a $\theta$-curve.
\end{Proof}

By using Lemma~\ref{GammaM+1} and 
lemmata in this section,
we obtain the following corollary:

\begin{corollary}
\LABEL{OvalTypeBTypeCTypeGTypeH}
If there exists a minimal chart $\Gamma$ of 
type $(m;2,3,2)$,
then each of $\Gamma_{m+1}$ and $\Gamma_{m+2}$
contains either the union of 
an oval and a skew $\theta$-curve, or
 one of four graphs as shown 
in  Fig.~\ref{Fig04}$($b$)$,$($c$)$,$($g$)$,$($h$)$.
\end{corollary}

%%%%%%%%%%%%%%%%%%%%%%%%
%%%%%%%%%%%%%%%%%%%%%%%%
%%%%%%%%%%%%%%%%%%%%%%%%
%\newpage

\section{RO-families of pseudo charts}
\LABEL{s:ROfamilies}

In this section,
we consider a minimal chart $\Gamma$ of type $(m;2,3,2)$
such that $\Gamma_{m+1}$ contains either an oval, or 
one of the four graphs as shown in
Fig.~\ref{Fig04}(b),(c),(g),(h).
We investigate 
that the chart $\Gamma$ contains what kind of pseudo charts.

Let $\Gamma$ be a chart, 
$D$ a disk, and 
$G$ a pseudo chart with $G \subset D$.
Let $r:D\to D$ be a reflection of $D$, and $G^*$ the pseudo chart obtained from $G$ by changing the orientations of all of the edges.
Then the set $\{G,G^*, r(G), r(G^*)\}$ 
is called the {\it RO-family of the pseudo chart $G$}.

In our argument,
we often need a name for an unnamed edge by using a given edge and a given white vertex.
For the convenience,
we use the following naming:
Let $e',e_i,e''$ be three consecutive edges containing  a white vertex $w_j$. Here, 
the two edges $e'$ and $e''$ are unnamed edges. 
There are six arcs in a neighborhood $U$ of the white vertex $w_j$. 
If the three arcs $e'\cap U$, $e_i \cap U$, $e'' \cap U$ lie anticlockwise around the white vertex $w_j$ in this order, 
then $e'$ and $e''$ are denoted by $a_{ij}$ and $b_{ij}$ 
respectively (see Fig.~\ref{Fig21}).
There is a possibility $a_{ij}=b_{ij}$ if they are contained in a loop.

%%%%%%%%%%%%%%%%%%
%%%%%%%%%%%%%%%%%% Figure
%%%%%%%%%%%%%%%%%%
\begin{figure}[htb]
\centerline{\includegraphics{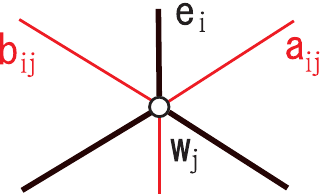}}
\caption{\LABEL{Fig21}
The three edges $a_{ij},e_i,b_{ij}$ are consecutive edges around the white vertex $w_j$.}
\end{figure}

%%%%%%%%%%%%%%%%%%%%%%%%%
%%%%%%%%%%%%%%%%%%%%%%%%%

\begin{lemma}
\LABEL{LemmaTypeC}
Let $\Gamma$ be a minimal chart of type $(m;2,3,2)$.
If $\Gamma_{m+1}$
contains the graph as shown in 
Fig.~\ref{Fig04}$($c$)$,
then $\Gamma$ contains one of the RO-family of
the pseudo chart as shown
in Fig.~\ref{Fig22}$($a$)$.
\end{lemma}

%%%%%%%%%%%%%%%%%%
%%%%%%%%%%%%%%%%%% Figure
%%%%%%%%%%%%%%%%%%
\begin{figure}[htb]
\centerline{\includegraphics{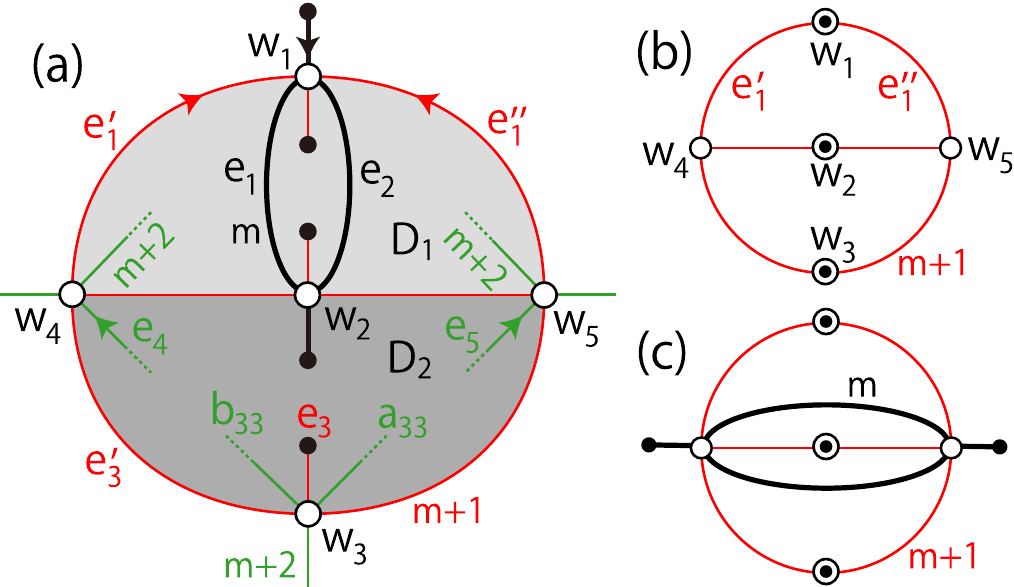}}
\caption{\LABEL{Fig22}
(a) A pseudo chart containing the graph as shown
in Fig.~\ref{Fig04}(c). 
The light gray region is $D_1$, and 
the dark gray region is $D_2$.
(b) $w_1,w_2,\cdots,w_5$ are white vertices.
(c) $w_4,w_5\in G$.
}
\end{figure}
\begin{Proof}
Suppose that $\Gamma_{m+1}$
contains the graph as shown
in Fig.~\ref{Fig04}$($c$)$.
We use the notations as shown in Fig.~\ref{Fig22}(b)
where $w_1,w_2,w_3$ are BW-vertices.
By Corollary~\ref{NoThetaGammaM},
the graph $\Gamma_m$ contains an oval $G$.
Thus by Lemma~\ref{LemmaBWvertexOval}(a),
there are two cases:
(i) $w_4,w_5\in G$ 
(see Fig.~\ref{Fig22}(c)),
(ii) the oval $G$ contains two of
BW-vertices $w_1,w_2,w_3$.

{\bf Case (i).}
By Lemma~\ref{LemmaOvalNoGraph}(a),
the chart $\Gamma$ is not minimal.
This is a contradiction.
Hence Case (i) does not occur.

{\bf Case (ii).}
Without loss of generality
we can assume $w_1,w_2\in G$.
Let $e_3$ be the terminal edge of label $m+1$ at $w_3$.
Let $e_1',e_1''$ be the internal edges of label $m+1$
at $w_1$ containing $w_4,w_5$ respectively
(see Fig.~\ref{Fig22}(b)).

Now, the graph in $\Gamma_{m+1}$ as shown in 
Fig.~\ref{Fig04}(c)
separates $S^2$ into three disks.
One of the three disks contains 
internal edges of label $m$ in $G$,
say $D_1$.
One of the three disks contains 
the terminal edge $e_3$,
say $D_2$.
The last one is denoted by $D_3$.
By Assumption~\ref{Infinity},
we can assume that 
the disk $D_3$ contains the point at infinity $\infty$.

If necessary
we change the orientation of all the edges of $\Gamma$,
we can assume that 
the terminal edge of label $m$ at $w_1$
is oriented inward at $w_1$.
Then by Assumption~\ref{AssumeTerminal}
\begin{enumerate}
\item[(1)]
the two  edges $e_1',e_1''$
 are oriented inward at $w_1$.
\end{enumerate}

Since $w_1,w_2\in \Gamma_m\cap\Gamma_{m+1}$, 
$w_3,w_4,w_5\in \Gamma_{m+1}$ and since
$\Gamma$ is of type $(m;2,3,2)$,
we have $w_3,w_4,w_5\in\Gamma_{m+1}\cap\Gamma_{m+2}$.
Let $e_4,e_5$ be the internal edges
(possibly terminal edges) of label $m+2$ at $w_4,w_5$
in $D_2$ respectively.
Then by (1),
the two edges $e_4,e_5$ 
are oriented inward at $w_4,w_5$
respectively.
Thus $\Gamma$ contains the pseudo chart as shown in 
Fig.~\ref{Fig22}(a).
\end{Proof}

\begin{lemma}
\LABEL{LemmaTypeB}
Let $\Gamma$ be a minimal chart of type $(m;2,3,2)$.
If $\Gamma_{m+1}$
contains the graph as shown in 
Fig.~\ref{Fig04}$($b$)$,
then $\Gamma$ contains one of the RO-families of
the two pseudo charts as shown
in Fig.~\ref{Fig23}.
\end{lemma}

%%%%%%%%%%%%%%%%%%
%%%%%%%%%%%%%%%%%% Figure
%%%%%%%%%%%%%%%%%%
\begin{figure}[htb]
\centerline{\includegraphics{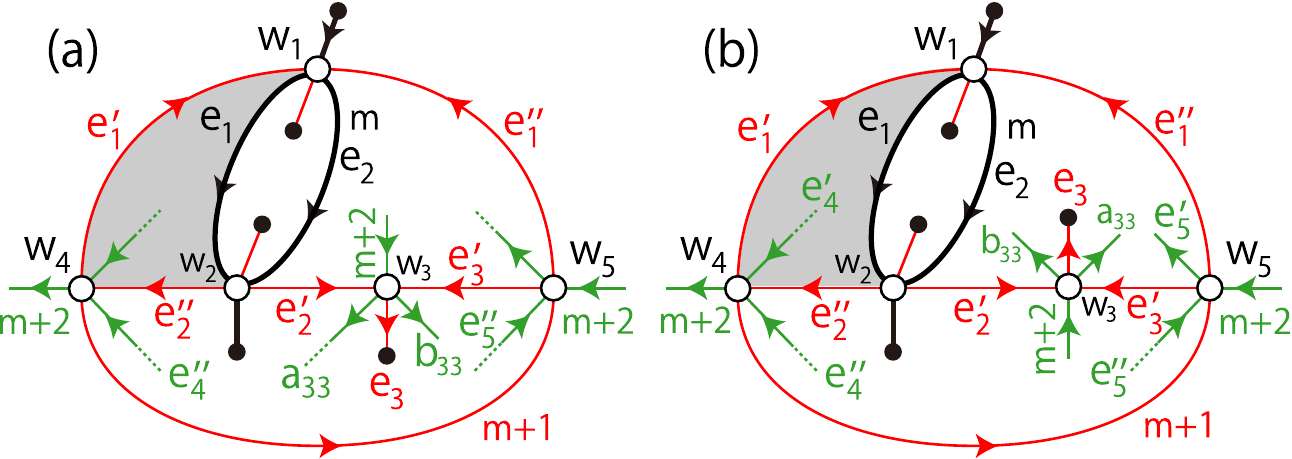}}
\caption{\LABEL{Fig23}
Pseudo charts containing the graph as shown
in Fig.~\ref{Fig04}(b).
The gray regions are the disk $D_1'$.
}
\end{figure}

\begin{Proof}
Suppose that $\Gamma_{m+1}$
contains the graph as shown
in Fig.~\ref{Fig04}$($b$)$.
We use the notations as shown in Fig.~\ref{Fig24}(a)
where $w_1,w_2,w_3$ are BW-vertices.
By Corollary~\ref{NoThetaGammaM},
the graph $\Gamma_m$ contains an oval $G$.
Thus by Lemma~\ref{LemmaBWvertexOval}(a),
there are three cases:
(i) $w_4,w_5\in G$ 
(see Fig.~\ref{Fig24}(b),(c),(d)),
(ii) $w_1,w_2\in G$ or $w_1,w_3\in G$,
(see Fig.~\ref{Fig24}(e)),
(iii) $w_2,w_3\in G$
(see Fig.~\ref{Fig24}(f)).

{\bf Case (i).}
Let $e',e'',e'''$ be the internal edges of label $m+1$
at $w_4$ containing
$w_1,w_2,w_5$,
respectively
(see Fig.~\ref{Fig24}(a)).
Let $D$ be the 2-angled disk of $\Gamma_m$
not containing the terminal edge of label $m$ at $w_4$.
There are three cases:
$e'\subset D$ or
$e''\subset D$ or
$e'''\subset D$.

If $e'\subset D$ (see Fig.~\ref{Fig24}(b)),
then
by Lemma~\ref{LemmaOvalNoGraph}(a)
the chart $\Gamma$ is not minimal.
This is a contradiction.
If $e''\subset D$ (see Fig.~\ref{Fig24}(c)),
then by Lemma~\ref{LemmaOvalLens}
there exists a lens in $Cl(S^2-D)$.
This contradicts Lemma~\ref{NoLens}.
If $e'''\subset D$ (see Fig.~\ref{Fig24}(d)),
then by Lemma~\ref{LemmaOvalLens}
there exists a lens in $D$.
This contradicts Lemma~\ref{NoLens}.
Hence Case (i) does not occur.

{\bf Case (ii).}
Without loss of generality
we can assume that $w_1,w_2\in G$.
Let $e_1,e_2$ be the internal edges of label $m$ in $G$.

Now, the graph in $\Gamma_{m+1}$ as shown in
Fig.~\ref{Fig04}(b)
separates $S^2$ into three disks.
One of the three disks contains 
both of $e_1$ and $e_2$.

Moreover, since $w_1,w_2\in\Gamma_m\cap\Gamma_{m+1}$,
 $w_3,w_4,w_5\in\Gamma_{m+1}$
and since $\Gamma$ is of type $(m;2,3,2)$,
we have $w_3,w_4,w_5\in \Gamma_{m+1}\cap\Gamma_{m+2}$.
Hence the chart $\Gamma$ contains 
the pseudo chart as shown in
Fig.~\ref{Fig24}(e).
We use the notations as shown in 
Fig.~\ref{Fig24}(e),
where 
$e_1',e_1''$
are internal edges of label $m+1$
at $w_1$,
$e_2',e_2''$
are internal edges of label $m+1$
at $w_2$ containing $w_3,w_4$ respectively,
$e_3'$ is an internal edge of label $m+1$
connecting $w_3,w_5$.

Without loss of generality
we can assume that
the terminal edge of label $m$ at $w_1$
is oriented inward at $w_1$.
Thus by Assumption~\ref{AssumeTerminal},
the two edges $e_1',e_1''$ 
are oriented inward at $w_1$, and
the two edges $e_1$ and $e_2$ 
are oriented from $w_1$ to $w_2$.
Hence
the two edges $e_2',e_2''$
are oriented outward at $w_2$.
Thus the edge $e_2'$
is oriented from $w_2$ to the BW-vertex $w_3$.
Hence by Lemma~\ref{OriBWvertex}
the edge $e_3'$
is oriented from $w_5$ to $w_3$.
Moreover,
we have the orientation of other edges.
Thus $\Gamma$ contains one of 
the two pseudo charts as shown
in Fig.~\ref{Fig23}.

{\bf Case (iii).}
By Lemma~\ref{LemmaOvalTwoEdges},
the two internal edges $e_1,e_2$ of label $m$ in $G$
are contained in the closure of 
the same connected component $F$ of $S^2-\Gamma_{m+1}$
(see Fig.~\ref{Fig24}(f)).
Thus the curve $e_1\cup e_2$ bounds 
a 2-angled disk of $\Gamma_m$ in $Cl(F)$,
say $D$.
Hence $Cl(S^2-D)$ is also a 2-angled disk of $\Gamma_m$,
and 
by Lemma~\ref{LemmaOvalLens} 
the disk $Cl(S^2-D)$ contains a lens.
This contradicts Lemma~\ref{NoLens}.
Hence Case (iii) does not occur.

Therefore $\Gamma$ 
contains one of the RO-families of 
 the two pseudo charts as shown in 
Fig.~\ref{Fig23}.
\end{Proof}

%%%%%%%%%%%%%%%%%%
%%%%%%%%%%%%%%%%%% Figure
%%%%%%%%%%%%%%%%%%
\begin{figure}[htb]
\centerline{\includegraphics{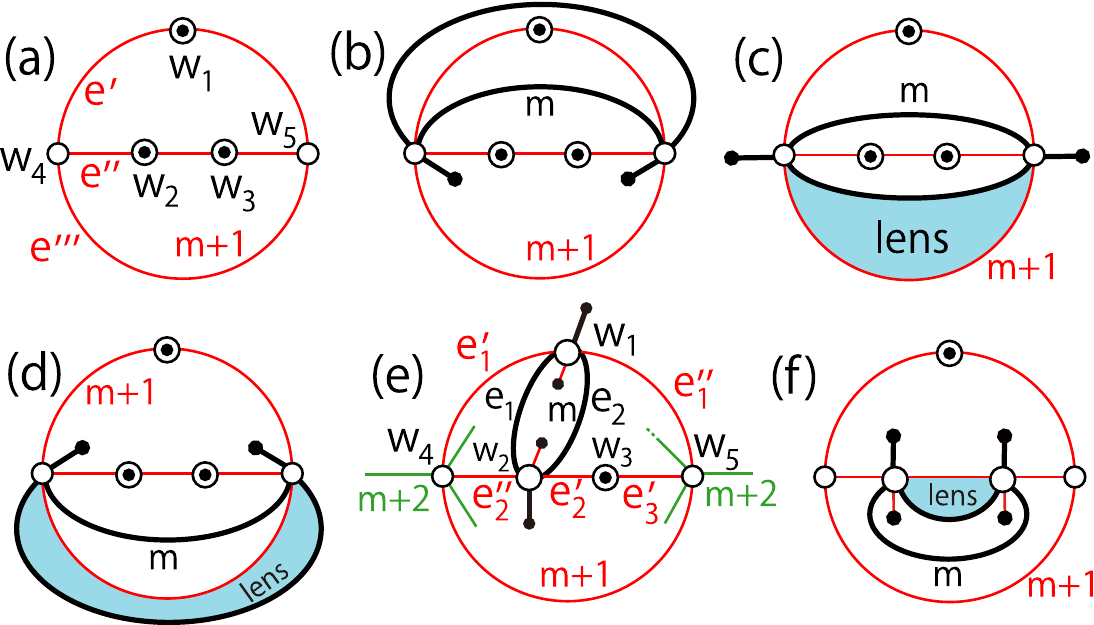}}
\caption{\LABEL{Fig24}
(a) $w_1,w_2,\cdots,w_5$ are white vertices.
(b),(c),(d) $w_4,w_5\in G$.
(e) $w_1,w_2\in G$.
(f) $w_2,w_3\in G$.
}
\end{figure}

The following lemma is not used in the this paper, but is used in the next paper \cite{ChartAppVII}.

\begin{lemma}
\LABEL{OvalGammaM+1Step0}
Let $\Gamma$ be a minimal chart 
of type $(m;2,3,2)$.
If $\Gamma_{m+1}$ contain an oval,
then 
$\Gamma$ contains one of the RO-families of
the two pseudo charts
as shown in Fig.~\ref{Fig25}$($a$)$ and $($b$)$.
\end{lemma}

\begin{figure}[htb]
\centerline{\includegraphics{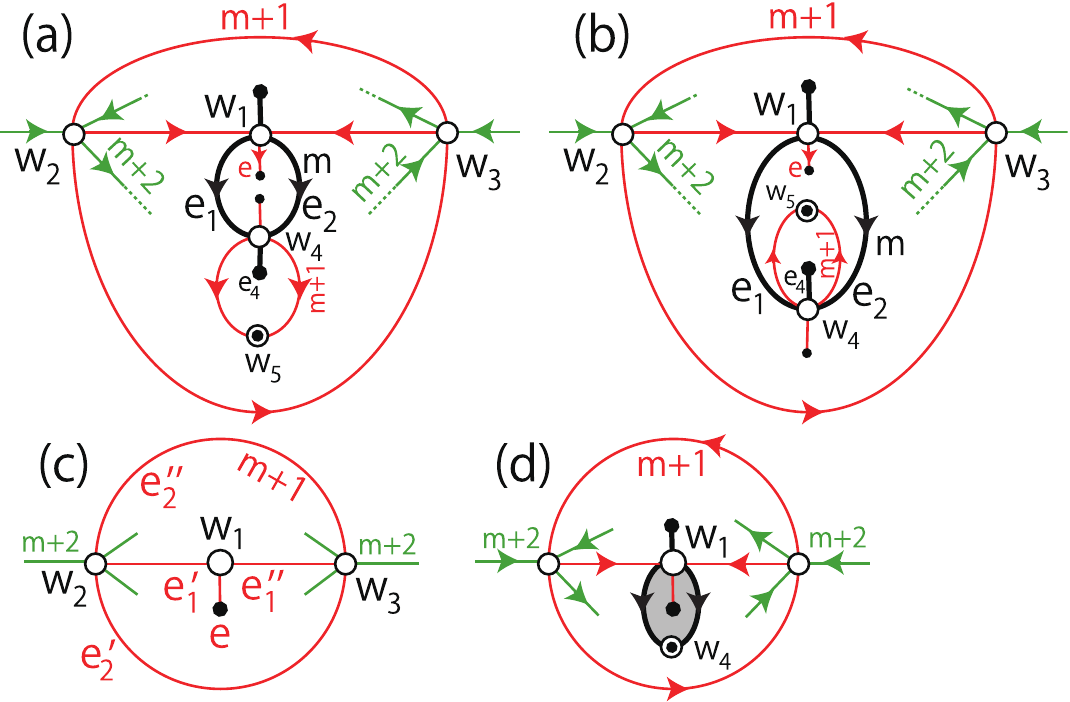}}
\caption{
\LABEL{Fig25}
(a),(b) Pseudo charts containing
a skew $\theta$-curve and an oval of label $m+1$.
(c) The skew $\theta$-curve of label $m+1$.
(d) The gray region is the disk $D$.
}
\end{figure}

\begin{Proof}
By Corollary~\ref{NoThetaGammaM},
the graph $\Gamma_m$ contains an oval $G$.
Since $\Gamma_{m+1}$ contain an oval,
by Lemma~\ref{GammaM+1}
the graph $\Gamma_{m+1}$ contains a skew $\theta$-curve.
Let $w_1,w_2,w_3$ be the white vertices in 
the skew $\theta$-curve such that
$w_1$ is a BW-vertex with respect to $\Gamma_{m+1}$.
Let $w_4,w_5$ be the white vertices in the oval 
of label $m+1$.
Then $w_4,w_5$ are BW-vertices 
with respect to $\Gamma_{m+1}$.
Hence by Lemma~\ref{LemmaBWvertexOval}(d),
there are three cases:
(i) $w_1,w_4\in G$ or $w_1,w_5\in G$,
(ii) $w_2,w_3\in G$,
(iii) $w_4,w_5\in G$.

{\bf Case (ii) and Case (iii).}
By Lemma~\ref{LemmaOval2angledDisk},
there exists a lens of $\Gamma$.
This contradicts Lemma~\ref{NoLens}.
Hence Case (ii) and Case (iii) do not occur.

{\bf Case (i).}
Without loss of generality
we can assume $w_1,w_4\in G$.
Since $w_1,w_4\in \Gamma_m\cap\Gamma_{m+1}$, 
$w_2,w_3\in \Gamma_{m+1}$ and since
$\Gamma$ is of type $(m;2,3,2)$,
we have $w_2,w_3\in\Gamma_{m+1}\cap\Gamma_{m+2}$
(see Fig.~\ref{Fig25}(c)).
We use the notations as shown 
in Fig.~\ref{Fig25}(c), where
$e$ is the terminal edge of label $m+1$ at $w_1$,
$e_1',e_1''$ are
internal edges of label $m+1$ at $w_1$,
and $e_2',e_2''$ are
internal edges of label $m+1$ 
connecting $w_2$ and $w_3$.

Without loss of generality,
we can assume that
\begin{enumerate}
\item[(1)] the terminal edge $e$ is 
oriented outward at $w_1$.
\end{enumerate}
Since the terminal edge $e$
is middle at $w_1$ by 
Assumption~\ref{AssumeTerminal},
the two edges $e_1',e_1''$ 
are oriented inward at $w_1$.
If necessary we reflect the chart $\Gamma$,
we can assume that
the edge $e_2'$ is
oriented from $w_2$ to $w_3$.
Looking at edges around $w_2$, 
the edge $e_2''$ is
oriented from $w_3$ to $w_2$.
Hence  
we have the orientation of the other edges 
of label $m+2$.
Let $e_1,e_2$ be the internal edges of label $m$ in $G$.
Then by (1)
the edges $e_1,e_2$ is oriented from $w_1$ to $w_4$.
Hence $\Gamma$ contains the pseudo chart
as shown in Fig.~\ref{Fig25}(d).

Let $D$ be the 2-angled disk of $\Gamma_m$
with $\partial D\ni w_1,w_4$
and $D\not\ni w_2$
(see Fig.~\ref{Fig25}(d)).
Let $e_4$ be the terminal edge of label $m$
at $w_4$.
There are two cases:
$e_4\not\subset D$ or $e_4\subset D$.
If $e_4\not\subset D$,
then the chart $\Gamma$ contains
the pseudo chart 
as shown in Fig.~\ref{Fig25}(a).
If $e_4\subset D$,
then the chart $\Gamma$ contains
the pseudo chart as shown in Fig.~\ref{Fig25}(b).
Therefore $\Gamma$ 
contains one of the RO-families of 
 the two pseudo charts as shown in 
Fig.~\ref{Fig25}(a),(b).
\end{Proof}

The following lemma is not used in the this paper, but is used in the next paper \cite{ChartAppVII}.

\begin{lemma}
\LABEL{LemmaTypeG}
Let $\Gamma$ be a minimal chart 
of type $(m;2,3,2)$.
If $\Gamma_{m+1}$ 
contains the graph as shown in 
Fig.~\ref{Fig04}$($g$)$,
then 
$\Gamma$ contains one of the RO-families of 
the two pseudo charts
as shown in Fig.~\ref{Fig26}$($a$)$,$($b$)$.
\end{lemma}

%%%%%%%%%%%%%%%%%%
%%%%%%%%%%%%%%%%%% Figure
%%%%%%%%%%%%%%%%%%
\begin{figure}[htb]
\centerline{\includegraphics{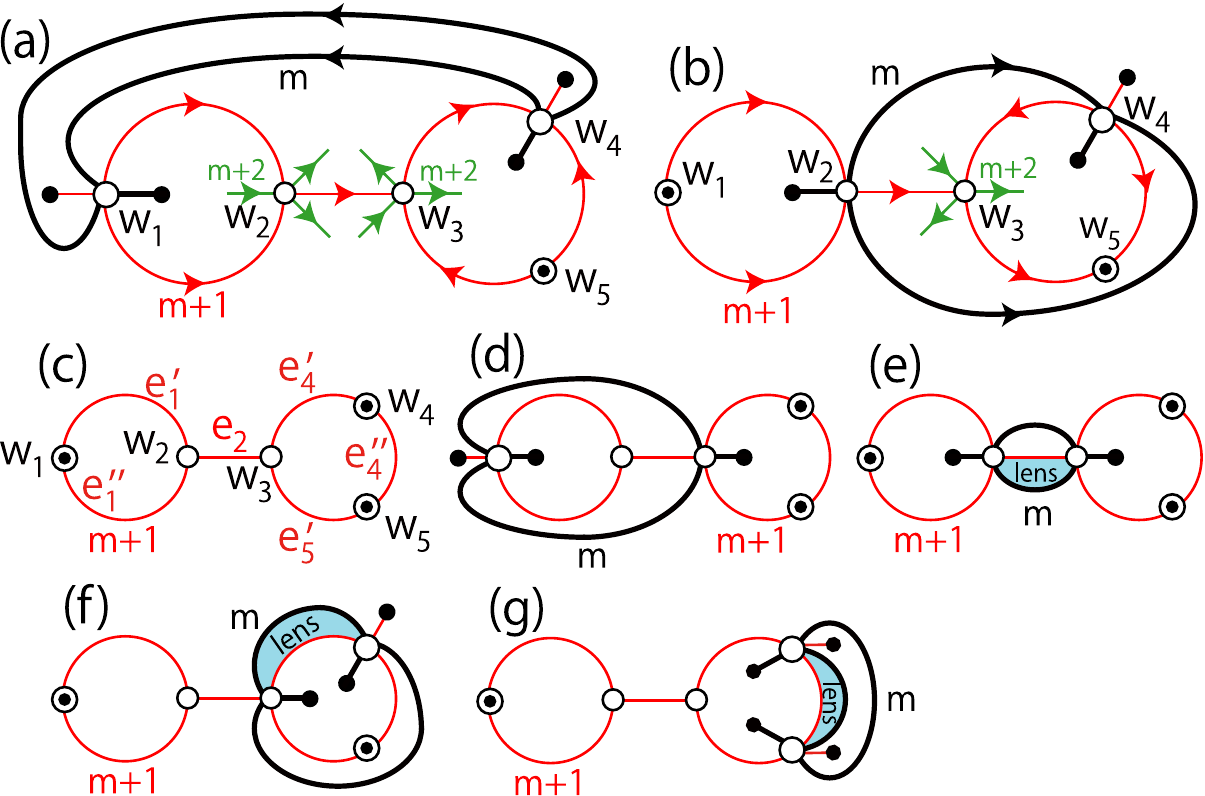}}
\caption{\LABEL{Fig26}
(a),(b) Pseudo charts containing the graph as shown in
Fig.~\ref{Fig04}(g).
(c) $w_1,w_2,\cdots,w_5$ are white vertices.
(d) $w_1,w_3\in G$.
(e) $w_2,w_3\in G$.
(f) $w_3,w_4\in G$.
(g) $w_4,w_5\in G$.}
\end{figure}

\begin{Proof}
Suppose that $\Gamma_{m+1}$
contains the graph as shown
in Fig.~\ref{Fig04}$($g$)$.
We use the notations as shown in Fig.~\ref{Fig26}(c)
where $w_1,w_4,w_5$ are BW-vertices.
By Corollary~\ref{NoThetaGammaM},
the graph $\Gamma_m$ contains an oval $G$.
There are seven cases:
(i) $w_1,w_2\in G$,
(ii) $w_1,w_3\in G$
(see Fig.~\ref{Fig26}(d)),
(iii) $w_1,w_4\in G$ or $w_1,w_5\in G$
(see Fig.~\ref{Fig26}(a)),
(iv) $w_2,w_3\in G$
(see Fig.~\ref{Fig26}(e)),
(v) $w_2,w_4\in G$ or $w_2,w_5\in G$
(see Fig.~\ref{Fig26}(b)),
(vi) $w_3,w_4\in G$ or $w_3,w_5\in G$
(see Fig.~\ref{Fig26}(f)),
(vii) $w_4,w_5\in G$
(see Fig.~\ref{Fig26}(g)).

{\bf Case (i).}
By Lemma~\ref{LemmaOval2angledDisk},
there exists a lens of $\Gamma$.
This contradicts Lemma~\ref{NoLens}.
Hence Case (i) does not occur.

{\bf Case (ii).}
By Lemma~\ref{LemmaOvalNoGraph}(b),
the chart $\Gamma$ is not minimal.
This is a contradiction.
Thus Case (ii) does not occur.

{\bf Case (iii).}
We use the notations as shown in Fig.~\ref{Fig26}(c)
where
$e_1',e_1''$ are two internal edges at $w_1$,
$e_2$ is the internal edge connecting
$w_2$ and $w_3$,
$e_4',e_4''$ are two internal edges at $w_4$,
and 
$e_5'$ is the internal edge 
connecting $w_3$ and $w_5$.

If necessary we reflect the chart,
we can assume that $w_1,w_4\in G$.
If necessary we change the orientation of all the edges,
we can assume that
the two edges $e_1',e_1''$
are oriented from $w_1$ to $w_2$.
Then
the edge $e_2$
is oriented from $w_2$ to $w_3$, and
the two internal edges of label $m$ in $G$
are oriented from $w_4$ to $w_1$.
Thus the two edges $e_4',e_4''$
are oriented inward at $w_4$.
Hence the edge $e_4''$ 
is oriented from $w_5$ to $w_4$.
Since $w_5$ is a BW-vertex with respect to $\Gamma_{m+1}$,
by Lemma~\ref{OriBWvertex}
the edge $e_5'$
is oriented from $w_5$ to $w_3$.
Moreover we have the orientation of the other edges.

Since $w_1,w_4\in\Gamma_m\cap\Gamma_{m+1}$,
$w_2,w_3\in \Gamma_{m+1}$
and since $\Gamma$ is of type $(m;2,3,2)$,
we have $w_2,w_3\in\Gamma_{m+1}\cap\Gamma_{m+2}$.
Therefore
$\Gamma$ contains the pseudo chart 
as shown in Fig.~\ref{Fig26}$($a$)$.

{\bf Case (iv).}
By Lemma~\ref{LemmaOvalLens},
there exists a lens.
This contradicts Lemma~\ref{NoLens}.
Thus Case (iv) does not occur.

{\bf Case (v).} 
If necessary we reflect the chart,
we can assume that $w_2,w_4\in G$. 
If necessary we change the orientation of all the edges,
we can assume that
the two internal edges $e_1',e_1''$ of label $m+1$ at $w_1$
are oriented from $w_1$ to $w_2$.
By a similar way as Case (iii),
we can show that
$\Gamma$ contains the pseudo chart
as shown in Fig.~\ref{Fig26}(b).

{\bf Case (vi).}
By Lemma~\ref{LemmaOvalTwoEdges},
the two internal edges $e_1,e_2$ of label $m$ in $G$
are contained in the closure of 
the same connected component $F$ of $S^2-\Gamma_{m+1}$.
Thus the curve $e_1\cup e_2$ bounds 
a 2-angled disk $D$ of $\Gamma_m$ with $w_5\in D$
(see Fig.~\ref{Fig26}(f)),
and the disk $D$ contains a lens 
by Lemma~\ref{LemmaOvalLens}.
This contradicts Lemma~\ref{NoLens}.
Thus Case (vi) does not occur.

{\bf Case (vii).}
By the similar way of Case (vi),
there exists a lens.
This contradicts Lemma~\ref{NoLens}.
Thus Case (vii) does not occur.

Therefore
$\Gamma$ contains one of the RO-families of
the two pseudo charts as shown in 
Fig.~\ref{Fig26}(a),(b).
\end{Proof}

The following lemma is not used in the this paper, but is used in the next paper \cite{ChartAppVII}.

\begin{lemma}
\LABEL{LemmaTypeH}
Let $\Gamma$ be a minimal chart 
of type $(m;2,3,2)$.
If $\Gamma_{m+1}$ contains the graph as shown in 
Fig.~\ref{Fig04}$($h$)$,
then 
$\Gamma$
contains one of the RO-families of the two pseudo charts
as shown in Fig.~\ref{Fig27}$($a$)$, $($b$)$.
\end{lemma}

%%%%%%%%%%%%%%%%%%
%%%%%%%%%%%%%%%%%% Figure
%%%%%%%%%%%%%%%%%%
\begin{figure}[htb]
\centerline{\includegraphics{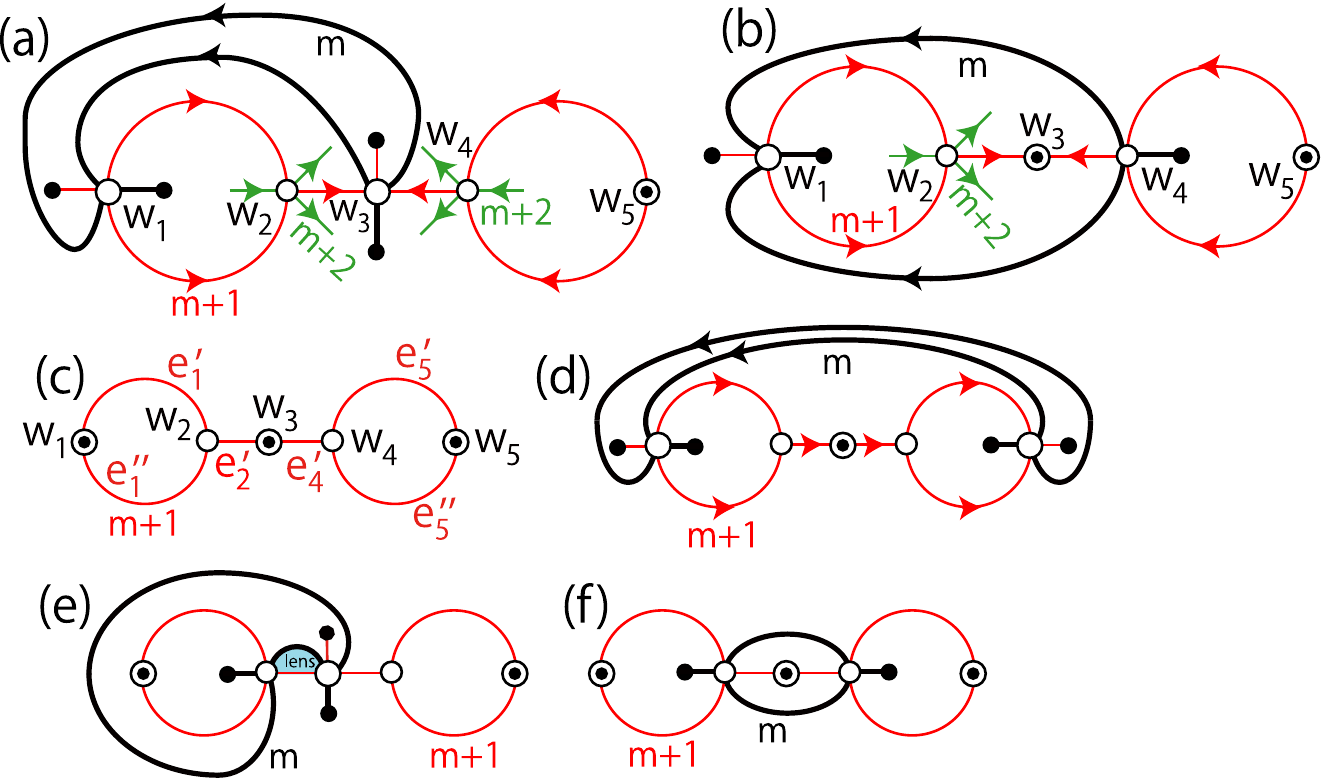}}
\caption{\LABEL{Fig27}
(a),(b) Pseudo charts containing the graph as
shown in Fig.~\ref{Fig04}(h). 
(c) $w_1,w_2,\cdots,w_5$ are white vertices.
(d) $w_1,w_5\in G$.
(e) $w_2,w_3\in G$.
(f) $w_2,w_4\in G$.}
\end{figure}

\begin{Proof}
Suppose that $\Gamma_{m+1}$
contains the graph as shown
in Fig.~\ref{Fig04}$($h$)$.
We use the notations as shown in 
Fig.~\ref{Fig27}(c)
where $w_1,w_3,w_5$ are BW-vertices,
$e_1',e_1''$ are two internal edges at $w_1$,
$e_2'$ is the internal edge connecting
$w_2$ and $w_3$,
$e_4'$ is the internal edge connecting
$w_3$ and $w_4$,
$e_5',e_5''$ are two internal edges at $w_5$.

By Corollary~\ref{NoThetaGammaM},
the graph $\Gamma_m$ contains an oval $G$.
There are six cases:
(i) $w_1,w_2\in G$ or $w_5,w_4\in G$,
(ii) $w_1,w_3\in G$ or $w_5,w_3\in G$
(see Fig.~\ref{Fig27}(a)),
(iii) $w_1,w_4\in G$ or $w_5,w_2\in G$
(see Fig.~\ref{Fig27}(b)),
(iv) $w_1,w_5\in G$
(see Fig.~\ref{Fig27}(d)),
(v) $w_2,w_3\in G$ or $w_4,w_3\in G$
(see Fig.~\ref{Fig27}(e)),
(vi) $w_2,w_4\in G$
(see Fig.~\ref{Fig27}(f)).

{\bf Case (i).}
By Lemma~\ref{LemmaOval2angledDisk},
there exists a lens of $\Gamma$.
This contradicts Lemma~\ref{NoLens}.
Hence Case (i) does not occur.

{\bf Case (ii).}
If necessary we reflect the chart,
we can assume that $w_1,w_3\in G$.
By Assumption~\ref{AssumeTerminal},
a neighborhood of $w_3$
contains the pseudo chart as shown in
Fig.~\ref{Fig14}.

If necessary we change the orientation of all the edges,
we can assume that
the two edges $e_1',e_1''$
are oriented from $w_1$ to $w_2$.
Then
the two internal edges of label $m$ in $G$
are oriented from $w_3$ to $w_1$, and
the two edges $e_2',e_4'$
are oriented inward at $w_3$.
Thus the edge $e_4'$ 
is oriented outward at $w_4$.
Hence by Lemma~\ref{OriBWvertex},
the two edges $e_5',e_5''$
are oriented from $w_5$ to $w_4$.
Thus we have the orientation of the other edges.

Since $w_1,w_3\in\Gamma_m\cap\Gamma_{m+1}$,
$w_2,w_4\in\Gamma_{m+1}$
and since $\Gamma$ is of type $(m;2,3,2)$,
we have $w_2,w_4\in\Gamma_{m+1}\cap\Gamma_{m+2}$.
Therefore
$\Gamma$ contains the pseudo chart
as shown in Fig.~\ref{Fig27}(a).

{\bf Case (iii).}
If necessary we reflect the chart,
we can assume that $w_1,w_4\in G$. 
If necessary we change the orientation of all the edges,
we can assume that
the two edges $e_1',e_1''$
are oriented from $w_1$ to $w_2$.
By a similar way as Case (ii),
we can show that
$\Gamma$ contains the pseudo chart
as shown in Fig.~\ref{Fig27}(b).

{\bf Case (iv).}
If necessary we change the orientation of all the edges,
we can assume that
the two edges $e_1',e_1''$
are oriented from $w_1$ to $w_2$.
Then 
\begin{enumerate}
\item[(1)] the edge $e_2'$
is oriented from $w_2$ to $w_3$
(i.e. the edge $e_2'$ is oriented inward at $w_3$), and 
\end{enumerate}
the two internal edges of label $m$ in $G$
are oriented from $w_5$ to $w_1$.
Thus
the two edges $e_5',e_5''$
are oriented from $w_4$ to $w_5$.
Hence the edge $e_4'$
is oriented from $w_3$ to $w_4$.
Thus the edge $e_4'$
is oriented outward at $w_3$ 
(see Fig.~\ref{Fig27}(d)).
However by (1) and Lemma~\ref{OriBWvertex}
we have a contradiction,
because $w_3$ is a BW-vertex 
with respect to $\Gamma_{m+1}$.
Thus Case (iv) does not occur.

{\bf Case (v).} 
Without loss of generality
we can assume $w_2,w_3\in G$.
If necessary we reflect the chart,
by Assumption~\ref{AssumeTerminal} 
(see Fig.~\ref{Fig14})
the chart $\Gamma$ contains the pseudo chart
as shown in Fig.~\ref{Fig27}(e).
Hence by Lemma~\ref{LemmaOvalLens},
there exists a lens.
This contradicts Lemma~\ref{NoLens}.
Thus Case (v) does not occur.

{\bf Case (vi).}
By Lemma~\ref{LemmaOvalNoGraph}(a),
the chart $\Gamma$ is not minimal.
This is a contradiction.
Thus Case (vi) does not occur.

Therefore
$\Gamma$ contains one of the RO-families of
the two pseudo charts as shown in 
Fig.~\ref{Fig27}(a),(b).
\end{Proof}

%%%%%%%%%%%%%%%%%%%%%%%%
%%%%%%%%%%%%%%%%%%%%%%%%
%%%%%%%%%%%%%%%%%%%%%%%%
%\newpage

\section{IO-Calculation}
\LABEL{s:IOC}

In this section,
we shall show that
neither $\Gamma_{m+1}$ nor $\Gamma_{m+2}$
contains the graph as shown in Fig.~\ref{Fig04}(c)
for any minimal chart $\Gamma$ of type $(m;2,3,2)$.

Let $\Gamma$ be a chart,
 and $v$ a vertex. 
Let $\alpha$ be a short arc of $\Gamma$ in a small neighborhood of $v$ such that $v$ is an endpoint of $\alpha$. 
If the arc $\alpha$ is oriented to $v$, then $\alpha$ is called {\it an inward arc}, 
and otherwise $\alpha$ is called {\it an outward arc}.

Let $\Gamma$ be an $n$-chart. 
Let $F$ be a closed domain with $\partial F\subset \Gamma_{k-1}\cup\Gamma_{k}\cup \Gamma_{k+1}$ for some label $k$ of $\Gamma$, where $\Gamma_0=\emptyset$ and $\Gamma_{n}=\emptyset$. 
By Condition (iii) for charts,
in a small neighborhood of each white vertex, there are three inward arcs and three outward arcs.
Also in a small neighborhood of each black vertex, there exists only one inward arc or one outward arc.
We often use the following fact, 
when we fix (inward or outward) arcs 
near white vertices and black vertices: 
\begin{enumerate}
\item[$(*)$]
{\it The number of inward arcs contained in $F\cap \Gamma_k$ is equal to the number of outward arcs in $F\cap \Gamma_k$.
}
\end{enumerate}
When we use this fact, 
we say that we use {\it IO-Calculation with respect to $\Gamma_k$ in $F$}.
For example, in a minimal chart $\Gamma$, 
consider the pseudo chart as shown in Fig.~\ref{Fig28} 
where
\begin{enumerate}
\item[(1)] $F$ is a $4$-angled disk of $\Gamma_{k-1}$,
\item[(2)]  $v_1,v_2,v_3,v_4$ are white vertices 
in $\partial F$  with  
$v_1\in\Gamma_{k-2}\cap\Gamma_{k-1}$ and
$v_2,v_3,v_4\in\Gamma_{k-1}\cap\Gamma_{k}$,
\item[(3)] $e_1$ is a terminal edge of label $k-2$ at $v_1$,
\item[(4)] $e_3$ is a terminal edge of label $k-1$ 
oriented inward at $v_3$,
\item[(5)]  for $i=2,4$,
the edge $e_i$ of label $k$ is oriented inward at $v_i$.
\end{enumerate}
Then we can show that $w(\Gamma\cap{\rm Int}F)\ge1$.
Suppose $w(\Gamma\cap{\rm Int}F)=0$.
By Assumption~\ref{AssumeTerminal},
the terminal edge $e_3$ contains a middle arc.
Thus
\begin{enumerate}
\item[(6)]  neither of edges $a_{33},b_{33}$ of label $k$ 
is middle at $v_3$
 (by Assumption~\ref{AssumeTerminal},
neither of them is a terminal edge). 
\end{enumerate}
Hence by (4),
\begin{enumerate}
\item[(7)] both of edges $a_{33},b_{33}$ of label $k$ 
are oriented inward at $v_3$. 
\end{enumerate}
If both of $e_2$ and $e_4$ are terminal edges of label $k$,
then 
by (5), (6), (7)
the number of inward arcs in $F\cap \Gamma_{k}$ is four,  
but the number of outward arcs in $F\cap \Gamma_{k}$ is two. 
This contradicts the fact $(*)$. 
Similarly if one of $e_2$ and $e_4$ is 
not a terminal edge of label $k$,
then we have the same contradiction.
Thus $w(\Gamma\cap{\rm Int}F)\ge1$.
Instead of the above argument, 
we just say that 
\begin{enumerate}
\item[]
{\it we have $w(\Gamma\cap{\rm Int}F)\ge1$ 
by IO-Calculation with respect to $\Gamma_{k}$ in $F$.}
\end{enumerate}

%%%%%%%%%%%%%%%%%%
%%%%%%%%%%%%%%%%%% Figure
%%%%%%%%%%%%%%%%%%
\begin{figure}
\centerline{\includegraphics{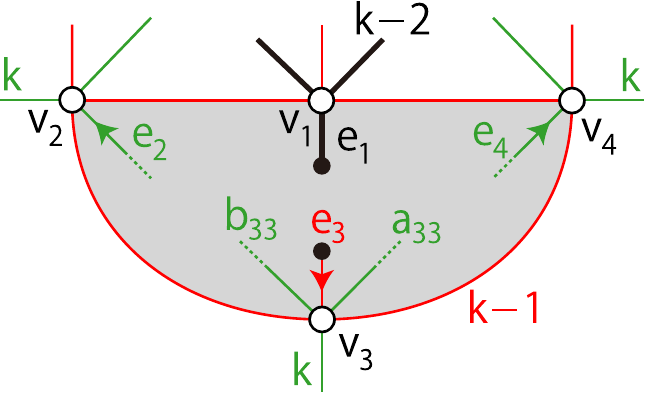}}
\caption{\LABEL{Fig28} The gray region is the disk $F$.}
\end{figure}

\begin{lemma}
{\rm $($\cite[Lemma 5.4]{ChartApp1}$)$} 
\LABEL{CorDiskLemma}
If a minimal chart $\Gamma$ contains 
the pseudo chart as shown in Fig.~\ref{Fig29}, 
then the interior of the disk $D^*$ 
contains at least one white vertex, where $D^*$ is the disk with the boundary $e_3^* \cup e_4^* \cup e^*$.
\end{lemma}

%%%%%%%%%%%%%%%%%%
%%%%%%%%%%%%%%%%%% Figure
%%%%%%%%%%%%%%%%%%
\begin{figure}[htb]
\centerline{\includegraphics{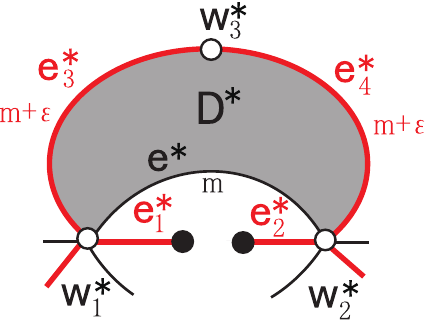}}
\caption{\LABEL{Fig29}
The gray region is the disk $D^*$.
 The label of the edge $e^*$ is $m$,
and $\varepsilon\in\{+1,-1\}$.}
\end{figure}

%%%%%%%%%%%%%%%%%%%%%%%%%
%%%%%%%%%%%%%%%%%%%%%%%%%

\begin{lemma}
\LABEL{NoMinimalTypeC}
Let $\Gamma$ be a minimal chart of type $(m;2,3,2)$.
Then neither $\Gamma_{m+1}$ nor $\Gamma_{m+2}$
contains the graph as shown in 
Fig.~\ref{Fig04}$($c$)$.
\end{lemma}

\begin{Proof}
Suppose that $\Gamma_{m+1}$ contains the graph as shown in 
Fig.~\ref{Fig04}$($c$)$.
By Lemma~\ref{LemmaTypeC},
the chart $\Gamma$ contains one of the RO-family of
the pseudo chart as shown in 
Fig.~\ref{Fig22}$($a$)$.
We use the notations as shown in
Fig.~\ref{Fig22}(a)
where
$e_1,e_2$ are the internal edges 
of label $m$ with $e_1\cap e_2\ni w_1,w_2$,
$e_3$ is the terminal edge of label $m+1$ at $w_3$,
and
\begin{enumerate}
\item[(1)] the two internal edges $e_4,e_5$ (possibly terminal edges) of label $m+2$ are oriented inward at
$w_4,w_5$, respectively.
\end{enumerate}

Now, the graph in $\Gamma_{m+1}$ as shown in 
Fig.~\ref{Fig04}(c)
separates $S^2$ into three disks.
One of them contains $e_1$,
say $D_1$,
and one of them contains 
the terminal edge $e_3$,
say $D_2$.

{\bf Claim 1.}
$w(\Gamma\cap{\rm Int}D_1)\ge2$ and
$w(\Gamma\cap{\rm Int}D_2)=0$.

{\it Proof of Claim~$1$.} 
The curve $e_1\cup e_2$ separates the disk $D_1$
into three disks.
One of them contains $w_4$, say $D_1'$,
and
one of them contains $w_5$, say $D_1''$.
Apply Lemma~\ref{CorDiskLemma}
considering as $D^*=D_1'$ and $w_3^*=w_4$,
we have 
$w(\Gamma\cap{\rm Int}D_1')\ge1$.
Similarly
we can show that 
$w(\Gamma\cap{\rm Int}D_1'')\ge1$.
Since $D_1\supset D_1'\cup D_1''$ and 
${\rm Int}D_1'\cap {\rm Int}D_1''=\emptyset$,
we have $w(\Gamma\cap{\rm Int}D_1)\ge2$.

Since $\Gamma$ is of type $(m;2,3,2)$,
we have $w(\Gamma)=7$ and $w(\Gamma_{m+1})=5$.
Thus
$$7=w(\Gamma)\ge w(\Gamma_{m+1})+w(\Gamma\cap{\rm Int}D_1)+w(\Gamma\cap{\rm Int}D_2)
\ge 5+2+w(\Gamma\cap{\rm Int}D_2).$$
Hence $w(\Gamma\cap{\rm Int}D_2)=0$.
Thus Claim~$1$ holds. {\hfill {$\square$}\vspace{1.5em}}

{\bf Claim 2.}
The terminal edge $e_3$
is oriented outward at $w_3$.

{\it Proof of Claim~$2$.} 
Suppose that 
$e_3$ is oriented inward at $w_3$.
Considering as $F=D_2$ and $k=m+2$
in the example of IO-Calculation
in Section~\ref{s:IOC},
the condition (1) implies that
we have $w(\Gamma\cap{\rm Int}D_2)\ge1$.
This contradicts the second equation of Claim~1.
Hence 
the terminal edge $e_3$ 
is oriented outward at $w_3$.
Thus Claim~2 holds. {\hfill {$\square$}\vspace{1.5em}}

Let $a_{33},b_{33}$ be the internal edges of label 
$m+2$ at $w_3$ in $D_2$
such that $a_{33},e_3,b_{33}$ lie anticlockwise 
around $w_3$ in this order
(see Fig.~\ref{Fig22}(a)).
By Assumption~\ref{AssumeTerminal},
\begin{enumerate}
\item[(2)] the terminal edge $e_3$ 
is middle at $w_3$.
\end{enumerate}

{\bf Claim 3.}
 $a_{33}=e_5$ and $b_{33}=e_4$.

{\it Proof of Claim~$3$.} 
By (2) and Assumption~\ref{AssumeTerminal},
 neither $a_{33}$ nor $b_{33}$ is
a terminal edge.
Moreover,
by Claim~2,
both of $a_{33}$ and $b_{33}$ 
are oriented outward at $w_3$. 
Thus by the second equation of Claim~$1$,
we have $a_{33}=e_5$ and $b_{33}=e_4$.
Hence Claim~3 holds. {\hfill {$\square$}\vspace{1.5em}}

Finally
we shall show that
there exists a lens of $\Gamma$.
Let $e_3'$ be the internal edge of label $m+1$
with $e_3'\ni w_3,w_4$
(see Fig.~\ref{Fig22}(a)).
By (2), neither $e_3'$ nor $b_{33}$ 
is middle at $w_3$.

By (2) and Claim~2,
the terminal edge $e_3$ is oriented outward at $w_3$
and middle at $w_3$.
Hence by Claim~3,
the edge $b_{33}=e_4$ is oriented from $w_3$ to $w_4$
and $e_3'$ is oriented from $w_4$ to $w_3$.
Thus  neither $e_3'$ nor $b_{33}$ 
is middle at $w_4$.
Hence
the curve $e_3'\cup b_{33}$ bounds a lens in $D_2$.
This contradicts Lemma~\ref{NoLens}.
Therefore
$\Gamma_{m+1}$ does not contain
the graph as shown in Fig.~\ref{Fig04}(c).

By Lemma~\ref{GammaM+1ToGammaM+2},
 we can show that
$\Gamma_{m+2}$ does not contain
the graph as shown in Fig.~\ref{Fig04}(c).
We complete the proof of Lemma~\ref{NoMinimalTypeC}.
\end{Proof}

%%%%%%%%%%%%%%%%%%%%%%%%

%\newpage

\section{Shifting Lemma}
\LABEL{s:ShiftingLemma}

In this section we shall show that
neither $\Gamma_{m+1}$ nor $\Gamma_{m+2}$
contains the graph as shown in Fig.~\ref{Fig04}(b)
for any minimal chart $\Gamma$ of type $(m;2,3,2)$.
Thus by Corollary~\ref{OvalTypeBTypeCTypeGTypeH} and
Lemma~\ref{NoMinimalTypeC},
we obtain the main theorem.

\begin{lemma}
\LABEL{NoMinimalTypeBTypeA}
Let $\Gamma$ be a chart of type $(m;2,3,2)$.
If $\Gamma$ contains the pseudo chart
as shown in Fig.~\ref{Fig23}$($a$)$,
then $\Gamma$ is not minimal.
\end{lemma}

\begin{Proof}
Suppose that $\Gamma$ is minimal.
We use the notations as shown in 
Fig.~\ref{Fig23}$($a$)$.
Here
$e_1,e_2$ are internal edges of label $m$, and
$e_4'',e_5'',a_{33},b_{33}$ are internal edges 
(possibly terminal edges) of label $m+2$
such that 
\begin{enumerate}
\item[(1)]
$e_4'',e_5''$ are oriented inward at $w_4,w_5$, 
respectively,
\item[(2)] 
$a_{33},b_{33}$ are oriented outward at $w_3$.
\end{enumerate}
Moreover
none of $e_4'',e_5'',a_{33},b_{33}$ are 
middle at $w_3,w_4$ or $w_5$.
Thus by Assumption~\ref{AssumeTerminal},
\begin{enumerate}
\item[(3)] 
none of $e_4'',e_5'',a_{33},b_{33}$ are terminal edges.
\end{enumerate}

Now, the graph $\Gamma_{m+1}$ 
contains the graph as shown in 
Fig.~\ref{Fig04}(b).
The graph in $\Gamma_{m+1}$ separates 
$S^2$ into three disks.
One of them contains 
the edges $e_1$ and $e_2$ of label $m$,
say $D_1$,
and one of them contains the edge $e_4''$,
say $D_2$.
Moreover,
the curve $e_1\cup e_2$ separates the disk $D_1$
into three disks.
One of them contains $w_4$,
say $D_1'$.
Apply Lemma~\ref{CorDiskLemma}
considering as $D^*=D_1'$ and $w_3^*=w_4$,
we have
\begin{enumerate}
\item[(4)] $w(\Gamma\cap{\rm Int}D_1')\ge1$.
\end{enumerate}

There are three cases:
(i) $w(\Gamma\cap {\rm Int}D_2)=0$,
(ii) $w(\Gamma\cap {\rm Int}D_2)=1$,
(iii) $w(\Gamma\cap {\rm Int}D_2)\ge2$.

{\bf Case (i).}
By using (1), (2) and (3),
we have $a_{33}=e_4''$ and $b_{33}=e_5''$.
Thus the curve $b_{33}\cup e_3'$ bounds a lens in $D_2$.
This contradicts Lemma~\ref{NoLens}.
Hence Case (i) does not occur.

{\bf Case (ii).}
Let $v$ be the white vertex in ${\rm Int}D_2$.
Since the five white vertices $w_1,w_2,\cdots,w_5$
are in $\Gamma_{m+1}$
and $\Gamma$ is of type $(m;2,3,2)$,
we have $v\in\Gamma_{m+2}\cap\Gamma_{m+3}$.
Thus
by using (1),(2) and (3),
we have a contradiction by IO-Calculation
with respect to $\Gamma_{m+2}$ in $D_2$.
Hence Case (ii) does not occur.

{\bf Case (iii).}
Since $\Gamma$ is of type $(m;2,3,2)$,
we have $w(\Gamma)=7$ and $w(\Gamma_{m+1})=5$.
Thus by (4) and 
the condition $w(\Gamma\cap {\rm Int}D_2)\ge2$ of this case,
$$7=w(\Gamma)\ge 
w(\Gamma_{m+1})+w(\Gamma\cap{\rm Int}D_1')
+w(\Gamma\cap{\rm Int}D_2)\ge 5+1+2=8.$$
This is a contradiction.
Hence Case (iii) does not occur.

Therefore the three cases do not occur.
Hence $\Gamma$ is not minimal.
\end{Proof}

Let $\Gamma$ and $\Gamma^\prime $ be C-move equivalent charts. 
Suppose that a pseudo chart $X$ of $\Gamma$ is also a pseudo chart of $\Gamma^\prime$. 
Then we say that 
$\Gamma$ is modified to $\Gamma^\prime$ by {\it C-moves keeping $X$ fixed}.
In Fig.~\ref{Fig30},
we give examples of C-moves keeping pseudo charts  fixed.

%%%%%%%%%%%%%%%%%%
%%%%%%%%%%%%%%%%%% Figure
%%%%%%%%%%%%%%%%%%
\begin{figure}[htb]
\centerline{\includegraphics{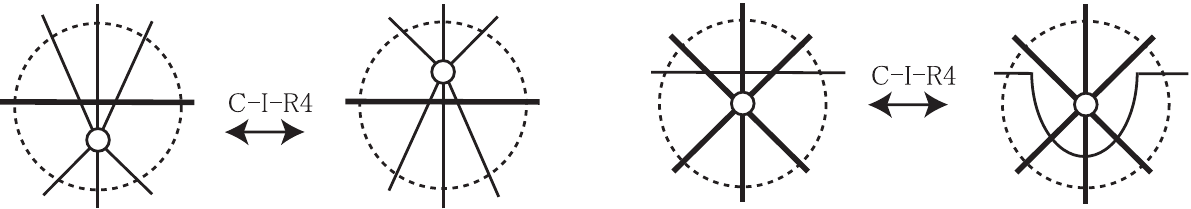}}
\caption{\LABEL{Fig30} 
C-moves keeping thicken figures fixed.}
\end{figure}

Let $\Gamma$ be a chart.
Let $\alpha$ be an arc in an edge of $\Gamma_m$, 
and 
$w$ a white vertex with $w\not\in\alpha$. 
Suppose that there exists an arc $\beta$ in $\Gamma$ 
such that
its end points are the white vertex $w$ 
and an interior point $p$ of the arc $\alpha$.
Then we say that 
{\it the white vertex $w$ connects with the point $p$ of $\alpha$ by the arc $\beta$}.

Let $\alpha$ be a simple arc,  
and $p,q$ points in $\alpha$. 
We denote by $\alpha[p,q]$ the subarc of $\alpha$ whose endpoints are $p$ and $q$.

\begin{lemma}
{\rm (\cite[Lemma 4.2]{ChartApp1})}
$($Shifting Lemma$)$
\LABEL{Shift} 
Let $\Gamma$ be a chart and
$\alpha$ an arc in an edge of $\Gamma_m$.
Let
$w$ be a white vertex of $\Gamma_k \cap\Gamma_{h}$ where $h=k+\varepsilon ,\varepsilon\in\{+1,-1\}$.
Suppose that 
the white vertex $w$ connects with a point $r$ 
of the arc $\alpha$ by an arc in an edge $e$ of $\Gamma_k$. 
Suppose that
one of the following two conditions is satisfied: 
\begin{enumerate}
\item[{\rm (1)}] $h>k>m$.
\item[{\rm (2)}] $h<k<m$.
\end{enumerate}
Then for any neighborhood $V$ of the arc $e[w,r]$
we can shift the white vertex $w$
to $e-e[w,r]$ 
along the edge $e$
by C-I-R2 moves,
C-I-R3 moves and C-I-R4 moves in $V$ 
keeping $\displaystyle 
\bigcup_{i<0}\Gamma_{k+i\varepsilon}
$ 
fixed $($see Fig.~\ref{Fig31}$)$.
\end{lemma}

%%%%%%%%%%%%%%%%%%
%%%%%%%%%%%%%%%%%% Figure
%%%%%%%%%%%%%%%%%%
\begin{figure}[htb]
\centerline{\includegraphics{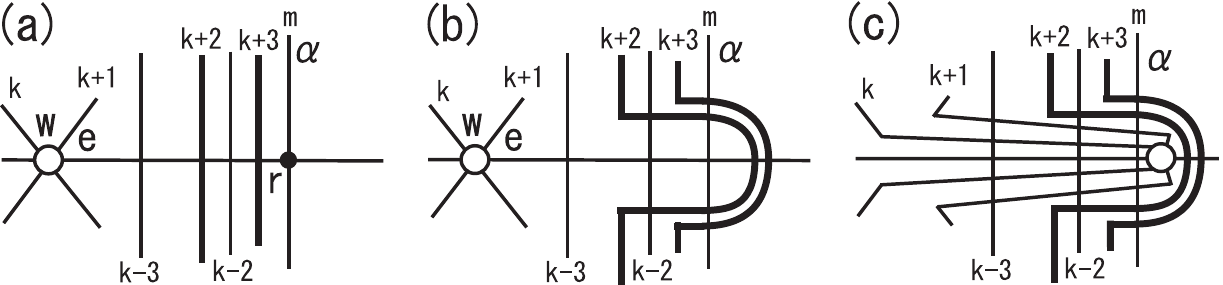}}
\caption{\LABEL{Fig31}
$k>m$ and $\varepsilon=+1$.}
\end{figure}

\begin{proposition}
\LABEL{LemmaNoMinimalTypeB}
Let $\Gamma$ be a minimal chart of type $(m;2,3,2)$.
Then neither $\Gamma_{m+1}$ nor $\Gamma_{m+2}$
contains the graph
as shown in Fig.~\ref{Fig04}$($b$)$.
\end{proposition}

%%%%%%%%%%%%%%%%%%%%%%

\begin{Proof}
Suppose that $\Gamma_{m+1}$ contains
 the graph
as shown in Fig.~\ref{Fig04}$($b$)$.
By Lemma~\ref{LemmaTypeB} and Lemma~\ref{NoMinimalTypeBTypeA},
the chart $\Gamma$ contains one of the RO-family
of the pseudo chart as shown in 
Fig.~\ref{Fig23}(b).
We use the notations as shown in 
Fig.~\ref{Fig23}(b).
Here,  $e_1,e_2$ are internal edges of label $m$, and
$e_4'',e_5',e_5'',a_{33},b_{33}$
are internal edges (possibly terminal edges) 
of label $m+2$ such that
\begin{enumerate}
\item[(1)] $e_4'',e_5''$ are oriented inward at $w_4,w_5$ respectively,
\item[(2)] $e_5',a_{33},b_{33}$ are oriented outward at $w_5,w_3,w_3$ respectively.
\end{enumerate}
Moreover,
none of $e_4'',e_5'',a_{33},b_{33}$
are middle at $w_3,w_4$ or $w_5$.
Thus by Assumption~\ref{AssumeTerminal},
\begin{enumerate}
\item[(3)] none of $e_4'',e_5'',a_{33},b_{33}$
are terminal edges.
\end{enumerate}

Now, the graph $\Gamma_{m+1}$ 
contains the graph as shown in 
Fig.~\ref{Fig04}(b).
The graph in $\Gamma_{m+1}$ separates 
$S^2$ into three disks.
One of them contains 
the edges $e_1$ and $e_2$,
say $D_1$,
and one of them contains the edge $e_4''$,
say $D_2$.
Moreover,
the curve $e_1\cup e_2$ separates the disk $D_1$
into three disks.
One of them contains $w_4$,
say $D_1'$.
Apply Lemma~\ref{CorDiskLemma}
considering as $D^*=D_1'$ and $w_3^*=w_4$,
we have
\begin{enumerate}
\item[(4)] $w(\Gamma\cap{\rm Int}D_1')\ge1$.
\end{enumerate}

{\bf Claim 1.} $w(\Gamma\cap {\rm Int}D_1)\ge1$ and 
$w(\Gamma\cap {\rm Int}D_2)\ge1$.

{\it Proof of Claim~$1$.} 
By (4) and $D_1'\subset D_1$,
we have $w(\Gamma\cap {\rm Int}D_1)\ge1$.

By (3),
neither $e_4''$ nor $e_5''$ is 
a terminal edge.
Since
$e_4'',e_5''$ are oriented inward at $w_4,w_5$ respectively
by (1), 
we have 
$w(\Gamma\cap{\rm Int}D_2)\ge1$
by IO-Calculation with respect to $\Gamma_{m+2}$
in $D_2$.
{\hfill {$\square$}\vspace{1.5em}}

{\bf Claim 2.}
$w(\Gamma\cap{\rm Int}D_1)=1$ and
$w(\Gamma\cap{\rm Int}D_1')=1$.

{\it Proof of Claim~$2$.} 
Suppose that $w(\Gamma\cap{\rm Int}D_1)\ge2$.
Since $\Gamma$ is of type $(m;2,3,2)$,
we have $w(\Gamma)=7$ and $w(\Gamma_{m+1})=5$.
Thus by the second inequality of Claim~1
$$7=w(\Gamma)\ge w(\Gamma_{m+1})+w(\Gamma\cap{\rm Int}D_1)+w(\Gamma\cap{\rm Int}D_2)\ge 5+2+1=8.$$
This is a contradiction.
Thus $w(\Gamma\cap{\rm Int}D_1)\le1$.
Hence by the first inequality of Claim~1,
we have $w(\Gamma\cap{\rm Int}D_1)=1$.

Thus by (4) and $D_1'\subset D_1$,
we have
$w(\Gamma\cap {\rm Int}D_1')=1$.
Hence Claim~2 holds. {\hfill {$\square$}\vspace{1.5em}}

Let $w_6$ be the white vertex in ${\rm Int}D_1'$.
Since the five white vertices $w_1,w_2,\cdots,w_5$
are in $\Gamma_{m+1}$ and
$\Gamma$ is of type $(m;2,3,2)$,
we have 
\begin{enumerate}
\item[(5)] $w_6\in\Gamma_{m+2}\cap\Gamma_{m+3}$.
\end{enumerate}

{\bf Claim~3.}
$a_{33}\ni w_6$ or $b_{33}\ni w_6$.

{\it Proof of Claim~$3$.} 
First we shall show that $a_{33}=e_4'$ or $a_{33}\ni w_6$.
By (3),
the edge $a_{33}$ is not a terminal edge.
Moreover,
by (2),
we have $a_{33}\not=e_5'$ and $a_{33}\not=b_{33}$.
Thus by Claim~2,
we have $a_{33}=e_4'$ or $a_{33}\ni w_6$.

Similarly we can show that
$b_{33}=e_4'$ or $b_{33}\ni w_6$.

If $a_{33}\not\ni w_6$ and $b_{33}\not\ni w_6$,
then $a_{33}=e_4'$ and $b_{33}=e_4'$.
This is a contradiction.
Therefore $a_{33}\ni w_6$ or $b_{33}\ni w_6$.
Thus Claim~3 holds. {\hfill {$\square$}\vspace{1.5em}}

If $a_{33}\ni w_6$,
let $e=a_{33}$,
otherwise let $e=b_{33}$.

{\bf Claim 4.}
We can move the white vertex $w_6$ 
from the disk $D_1'$ to the outside of $D_1'$.

{\it Proof of Claim~$4$.} 
Since the edge $e$ connects 
the vertex $w_6$ in ${\rm Int}D_1'$
and 
the vertex $w_3$ in the outside of $D_1'$,
the edge $e$ intersects the boundary $\partial D_1'$.
Let $x$ be the point in the edge $e$ 
with $e[w_6,x]\cap\partial D_1'=x$.
Since the edge $e$ is of 
label $m+2$
and since $\partial D_1'$ consists of 
the edge $e_1$ of label $m$ and
two internal edges of label $m+1$,
we have 
$$e[w_6,x]\cap e_1=e[w_6,x]\cap\partial D_1'=x.$$
Thus by (5),
the white vertex $w_6\in\Gamma_{m+2}\cap\Gamma_{m+3}$
connects with the point $x$ in the edge $e_1$ of label $m$
by the arc $e[w_6,x]$ of label $m+2$.
Hence by Shifting Lemma (Lemma~\ref{Shift}),
we can shift the white vertex $w_6$ by
C-I-R2 moves, C-I-R3 moves and C-I-R4 moves
in a neighborhood of the arc $e[w_6,x]$ 
keeping $\displaystyle{\bigcup_{i<0}\Gamma_{m+2+i}}$ fixed.
Thus we can shift the white vertex $w_6$ 
to the outside of $D_1'$ by C-moves
keeping $\partial D_1'$ fixed.
Therefore Claim~$4$ holds.
{\hfill {$\square$}\vspace{1.5em}}

By Claim~2 and Claim~4,
we have $w(\Gamma\cap{\rm Int}D_1')=0$.
However
we have a contradiction
by Lemma~\ref{CorDiskLemma}
considering as $D^*=D_1'$ and $w_3^*=w_4$.
Therefore $\Gamma_{m+1}$ does not contain
the graph as shown in Fig.~\ref{Fig04}(b).

By Lemma~\ref{GammaM+1ToGammaM+2},
 we can show that
$\Gamma_{m+2}$ does not contain
the graph as shown in Fig.~\ref{Fig04}(b).
We complete the proof of 
Proposition~\ref{LemmaNoMinimalTypeB}.
\end{Proof}

By Corollary~\ref{OvalTypeBTypeCTypeGTypeH},
Lemma~\ref{NoMinimalTypeC} and
Proposition~\ref{LemmaNoMinimalTypeB},
we have the main theorem (Theorem~\ref{OvalTypeGTypeH}).

%\newpage

%%%%%%%%%%%%%%%%%%%%%%%%
%%%%%%%%%%%%%%%%%%%%%%%%
%%%%%%%%%%%%%%%%%%%%%%%%

%%%%%%%%%%%%%%%%%%%%%%%%%%
%%%   Reference
%%%%%%%%%%%%%%%%%%%%%%%%%%

%%%%%%%%%%%%%%%%%%%%%%
%%%%%%%%%%%%%%%%%%%%%%%
%%%%%%%%%%%%%%%%%%%%%%%

\vspace{5mm}

\begin{minipage}{65mm}
{Teruo NAGASE
\\
{\small Tokai University \\
4-1-1 Kitakaname, Hiratuka \\
Kanagawa, 259-1292 Japan\\
\\
nagase@keyaki.cc.u-tokai.ac.jp
}}
\end{minipage}
\begin{minipage}{65mm}
{Akiko SHIMA 
\\
{\small Department of Mathematics, 
\\
Tokai University
\\
4-1-1 Kitakaname, Hiratuka \\
Kanagawa, 259-1292 Japan\\
shima@keyaki.cc.u-tokai.ac.jp
}}
\end{minipage}

\vspace{1cm}

{\bf List of terminologies}\vspace{2mm}\\
{\small $
\begin{array}{ll||}
\text{$k$-angled disk} & p11 \\
\text{BW-vertex} & p2 \\
\text{C-move~equivalent} & p5 \\
\text{chart} & p4 \\
\text{complexity} & p5 \\
\text{free edge} & p5 \\
\text{hoop} & p5 \\
\text{internal edge} & p9 \\
\text{inward} & p5 \\
\text{inward arc} & p31 \\
\text{IO-Calculation} & p31 \\
\text{keeping $X$ fixed} & p35 \\
\text{lens} & p10 \\
\text{loop} & p7\\
\text{middle arc} & p4 \\
\end{array}
~~
\begin{array}{ll}
\text{middle at $v$} & p5 \\
\text{minimal chart} & p5 \\
\text{outward} & p5 \\
\text{outward arc} & p31 \\
\text{oval} & p2 \\
\text{point at infinity $\infty$} & p5 \\
\text{pseudo chart} & p12 \\
\text{ring} & p5 \\
\text{RO-family} & p23 \\
\text{simple hoop} & p5 \\
\text{skew $\theta$-curve} & p2 \\
\text{terminal edge} & p2 \\
\text{type $(m;n_1,n_2,\cdots,n_k)$} & p1 \\
\text{$w$ connects with $p$ by an arc $\beta$} & p35 \\
\text{$\theta$-curve} & p2 \\
\end{array}
$}

\vspace{0.5cm}

{\bf List of notations}\vspace{2mm}\\
{\small $
\begin{array}{ll}
\text{$\Gamma_m$} & p1 \\
\text{$w(X)$} & p2 \\
\text{${\rm Int}X$} & p6 \\
\text{$\partial X$} & p6 \\
\text{$Cl(X)$} & p6 \\
\text{$a_{ij}$,$b_{ij}$} & p23 \\
\text{$\alpha[p,q]$} & p35 \\
\end{array}
$
}

\end{document}